\documentclass{amsart}
\makeatletter
\@ifundefined{lhs2tex.lhs2tex.sty.read}%
  {\@namedef{lhs2tex.lhs2tex.sty.read}{}%
   \newcommand\SkipToFmtEnd{}%
   \newcommand\EndFmtInput{}%
   \long\def\SkipToFmtEnd#1\EndFmtInput{}%
  }\SkipToFmtEnd

\newcommand\ReadOnlyOnce[1]{\@ifundefined{#1}{\@namedef{#1}{}}\SkipToFmtEnd}
\usepackage{amstext}
\usepackage{amssymb}
\usepackage{stmaryrd}
\DeclareFontFamily{OT1}{cmtex}{}
\DeclareFontShape{OT1}{cmtex}{m}{n}
  {<5><6><7><8>cmtex8
   <9>cmtex9
   <10><10.95><12><14.4><17.28><20.74><24.88>cmtex10}{}
\DeclareFontShape{OT1}{cmtex}{m}{it}
  {<-> ssub * cmtt/m/it}{}

\DeclareFontShape{OT1}{cmtt}{bx}{n}
  {<5><6><7><8>cmtt8
   <9>cmbtt9
   <10><10.95><12><14.4><17.28><20.74><24.88>cmbtt10}{}
\DeclareFontShape{OT1}{cmtex}{bx}{n}
  {<-> ssub * cmtt/bx/n}{}
\newcommand{\Conid}[1]{\mathit{#1}}
\newcommand{\Varid}[1]{\mathit{#1}}
\newcommand{\anonymous}{\kern0.06em \vbox{\hrule\@width.5em}}
\newcommand{\plus}{\mathbin{+\!\!\!+}}

\renewcommand{\leq}{\leqslant}
\renewcommand{\geq}{\geqslant}
\usepackage{polytable}

\@ifundefined{mathindent}%
  {\newdimen\mathindent\mathindent\leftmargini}%
  {}%

\def\resethooks{%
  \global\let\SaveRestoreHook\empty
  \global\let\ColumnHook\empty}
\newcommand*{\savecolumns}[1][default]%
  {\g@addto@macro\SaveRestoreHook{\savecolumns[#1]}}
\newcommand*{\restorecolumns}[1][default]%
  {\g@addto@macro\SaveRestoreHook{\restorecolumns[#1]}}
\newcommand*{\aligncolumn}[2]%
  {\g@addto@macro\ColumnHook{\column{#1}{#2}}}

\resethooks

\newcommand{\onelinecommentchars}{\quad-{}- }
\newcommand{\commentbeginchars}{\enskip\{-}
\newcommand{\commentendchars}{-\}\enskip}

\newcommand{\visiblecomments}{%
  \let\onelinecomment=\onelinecommentchars
  \let\commentbegin=\commentbeginchars
  \let\commentend=\commentendchars}

\newcommand{\invisiblecomments}{%
  \let\onelinecomment=\empty
  \let\commentbegin=\empty
  \let\commentend=\empty}

\visiblecomments

\newlength{\blanklineskip}
\newcommand{\hsindent}[1]{\quad}% default is fixed indentation

\makeatother
\EndFmtInput
\usepackage[latin1]{inputenc}
\usepackage{amsmath}
\usepackage{amsthm}
\usepackage{float}
\usepackage{url}

\newtheorem{theorem}{Theorem}[section]
\newtheorem{definition}[theorem]{Definition}

\floatstyle{ruled}
\newfloat{example}{ht!}{exa}
\floatname{example}{Example}

\begin{document}

\title{Enumerating the Saneblidze-Umble diagonal terms}

\author{Mikael Vejdemo-Johansson} 
\email{mik@math.uni-jena.de}
\address{Lehrstuhl für Algebra und Zahlentheorie \\
  Mathematisches Institut \\
  Fakultät für Mathematik und Informatik \\
  FSU Jena \\
  07737 Jena \\
  Germany} 

\thanks{The author acknowledges travel support from DFG Sachbehilfe
  grant GR 1585/4-1.}

%\classification{55-04 (Primary) 52B05 18D50}
%\keywords{permutahedron, associahedron, computer implementation, Saneblidze-Umble diagonal}

\begin{abstract}
  The author presents a computer implementation, calculating the terms
  of the Saneblidze-Umble diagonals on the permutahedron and the
  associahedron. The code is analyzed for correctness and presented in
  the paper, the source code of which simultaneously represents both
  the paper and the program.
\end{abstract}

%\received{~}
%\revised{~}
%\published{~}
%\submitted{~}

%\volumeyear{~}
%\volumenumber{~}
%\issuenumber{~}

%\startpage{1}

\maketitle

\section{Introduction}

In \cite{saneblidze-umble-2004}, Samson Saneblidze and Ron Umble gave
a combinatorial description of a diagonal on the permutahedron $P_n$,
together with the associated diagonal on the associahedron $K_n$
induced by the projection $P_n\to K_n$ due to Andy Tonks
\cite{tonks97}. This provides a tool for describing an
$A_\infty$-structure on the tensor product of two
$A_\infty$-(co)algebras.

There has been work done by Steve Weaver, in \cite{weaver05}, to
produce a computer implementation of the resulting algorithm for
listing the terms. However, this implementation has not been
widely disseminated, and further offers no transparency as to
how the code written corresponds to the mathematical
description. Furthermore, Andy Tonks \cite{tonks_personal07} has a
private implementation of the resulting algorithm.

Haskell is one of the stricter functional programming languages. It
also has very strong semantics. This means that a Haskell programmer
can reason about the written program in a way that would yield a
stringent proof of correctness.  In this paper, the author presents a
Haskell implementation of the enumeration algorithm of the
Saneblidze-Umble diagonal, interleaved with a description of the
diagonal itself, in such a manner that the resulting code is clear
enough and analogous enough to the mathematical argumentation for the
correctness of the code to be clear.

In section \ref{sec:sanebl-umble-diag}, the author presents the code
and the diagonal in parallel, beginning with the diagonal on the
permutahedron and finalizing with a discussion and implementation of
the projection. In section \ref{sec:inst-usage-example}, a practical
discussion of the usage of this paper for actual computation is given
as well as an explicit example session with one of the mainstream
Haskell interpreters. Finally, in section \ref{sec:haskell-notation},
an overview of the used symbols and library functions is given for
reference for the readers.

\section{The Saneblidze-Umble diagonal}
\label{sec:sanebl-umble-diag}

In this section we shall, in parallel, describe the diagonal on the
permutahedron and the resulting diagonal on the associahedron and give
computer code in the functional programming language Haskell to work
with computer calculations of the diagonal.

Since we shall be building a working programming library with code
performing the relevant calculations, we need to declare the library
as such.
\begingroup\par\noindent\advance\leftskip\mathindent\(
\begin{pboxed}\SaveRestoreHook
\column{B}{@{}l@{}}
\column{E}{@{}l@{}}
\>[B]{}\mathbf{module}\;\Conid{SaneblidzeUmbleSigns}\;\mathbf{where}{}\<[E]
\\
\>[B]{}\mathbf{import}\;\Conid{\Conid{Data}.List}{}\<[E]
\\
\>[B]{}\mathbf{import}\;\Conid{\Conid{Data}.Maybe}{}\<[E]
\\
\>[B]{}\mathbf{import}\;\Varid{qualified}\;\Conid{\Conid{Data}.Map}\;\Varid{as}\;\Conid{Map}{}\<[E]
\\
\>[B]{}\mathbf{import}\;\Conid{\Conid{Data}.Map}\;((\mathbin{!})){}\<[E]
\ColumnHook
\end{pboxed}
\)\par\noindent\endgroup\resethooks

We shall be working a lot with representations of faces of the
permutahedra $P_n$. These faces are indexed by ordered partitions, as
defined in the next section. More specifically, we shall be spending
most of our time working with $C_*(P_n) = C_*(P_n;k)$ -- the cellular
chain complex of the $n$th permutahedron with coefficients in some
field $k$. The aim of the efforts is to provide a cellular chain
homotopic to a diagonal of $P_n\times P_n$, which thus will be given
in $C_*(P_n\times P_n)=C_*(P_n)\otimes C_*(P_n)$.

\subsection{The diagonal on the permutahedron}
\label{sec:diag-perm}

The enumeration given by Saneblidze and Umble
\cite{saneblidze-umble-2004} is obtained by defining various matrices and
operations on those. The closure of the defined operations is in
bijective correspondence with the terms of a diagonal on the
permutahedron. This correspondence works by translating the matrices
into pairs of ordered partitions, indexing a pair of faces on the
permutahedron. 

In order to make the exposition here closer to the programming, I
shall omit further mention of matrices, and instead work exclusively
with the partitions. 

\begin{definition}
  An \emph{ordered partition} of the ordered set $\pi$ is a
  decomposition of $\pi$ into disjoint ordered subsets $\pi_1,\dots,\pi_n$
  covering $\pi$. Such a decomposition is denoted by
  $\pi_1|\pi_2|\dots|\pi_n$.

  Given a permutation $\sigma\in S_n$, the \emph{rising partition of
    $\sigma$} of $\sigma$ is a sequence of integer sequences, each of
  which is a maximal monotonically increasing subsequence of the
  $\sigma(1),\dots,\sigma(n)$ and the concatenation of all sequences
  form the sequence $\sigma(1),\dots,\sigma(n)$.

  Analogously, we define the \emph{falling partition of $\sigma$} of a
  permutation.
\end{definition}

Since the process of finding the rising partition and the falling
partition of $\sigma$ has a common abstraction, we shall here give the
code for that abstraction, and then the two relevant specializations.
\begingroup\par\noindent\advance\leftskip\mathindent\(
\begin{pboxed}\SaveRestoreHook
\column{B}{@{}l@{}}
\column{3}{@{}l@{}}
\column{5}{@{}l@{}}
\column{10}{@{}l@{}}
\column{11}{@{}l@{}}
\column{17}{@{}l@{}}
\column{34}{@{}l@{}}
\column{E}{@{}l@{}}
\>[B]{}\mathbf{type}\;\Conid{Sequence}{}\<[17]
\>[17]{}\mathrel{=}[\mskip1.5mu \Conid{Int}\mskip1.5mu]{}\<[E]
\\
\>[B]{}\mathbf{type}\;\Conid{Partition}{}\<[17]
\>[17]{}\mathrel{=}[\mskip1.5mu \Conid{Sequence}\mskip1.5mu]{}\<[E]
\\[\blanklineskip]
\>[B]{}\Varid{monotonicSequence}\mathbin{::}(\Varid{a}\to \Varid{a}\to \Conid{Bool})\to [\mskip1.5mu \Varid{a}\mskip1.5mu]\to [\mskip1.5mu [\mskip1.5mu \Varid{a}\mskip1.5mu]\mskip1.5mu]{}\<[E]
\\
\>[B]{}\Varid{monotonicSequence}\;\anonymous \;[\mskip1.5mu \mskip1.5mu]{}\<[34]
\>[34]{}\mathrel{=}[\mskip1.5mu \mskip1.5mu]{}\<[E]
\\
\>[B]{}\Varid{monotonicSequence}\;\anonymous \;[\mskip1.5mu \Varid{x}\mskip1.5mu]{}\<[34]
\>[34]{}\mathrel{=}[\mskip1.5mu [\mskip1.5mu \Varid{x}\mskip1.5mu]\mskip1.5mu]{}\<[E]
\\
\>[B]{}\Varid{monotonicSequence}\;\Varid{cmp}\;(\Varid{x}\mathbin{:}\Varid{y}\mathbin{:}\Varid{etc}){}\<[34]
\>[34]{}\mathrel{=}{}\<[E]
\\
\>[B]{}\hsindent{5}{}\<[5]
\>[5]{}\mathbf{if}\;{}\<[11]
\>[11]{}\Varid{x}\mathbin{`\Varid{cmp}`}\Varid{y}{}\<[E]
\\
\>[B]{}\hsindent{5}{}\<[5]
\>[5]{}\mathbf{then}\;{}\<[11]
\>[11]{}(\Varid{x}\mathbin{:}\Varid{s})\mathbin{:}\Varid{ss}{}\<[E]
\\
\>[B]{}\hsindent{5}{}\<[5]
\>[5]{}\mathbf{else}\;{}\<[11]
\>[11]{}[\mskip1.5mu \Varid{x}\mskip1.5mu]\mathbin{:}(\Varid{s}\mathbin{:}\Varid{ss}){}\<[E]
\\
\>[B]{}\hsindent{3}{}\<[3]
\>[3]{}\mathbf{where}{}\<[E]
\\
\>[3]{}\hsindent{2}{}\<[5]
\>[5]{}(\Varid{s}\mathbin{:}\Varid{ss})\mathrel{=}\Varid{monotonicSequence}\;\Varid{cmp}\;(\Varid{y}\mathbin{:}\Varid{etc}){}\<[E]
\\[\blanklineskip]
\>[B]{}\Varid{rising}\mathbin{::}\Conid{Sequence}\to \Conid{Partition}{}\<[E]
\\
\>[B]{}\Varid{rising}{}\<[10]
\>[10]{}\mathrel{=}\Varid{monotonicSequence}\;(\leq ){}\<[E]
\\[\blanklineskip]
\>[B]{}\Varid{falling}\mathbin{::}\Conid{Sequence}\to \Conid{Partition}{}\<[E]
\\
\>[B]{}\Varid{falling}{}\<[10]
\>[10]{}\mathrel{=}\Varid{monotonicSequence}\;(\geq ){}\<[E]
\ColumnHook
\end{pboxed}
\)\par\noindent\endgroup\resethooks

Now, a basis element of the tensor product $C_*(P_n)\otimes C_*(P_n)$
has the form $u\otimes v$ where $u$ and $v$ are both partitions of the
kind discussed here. Saneblidze and Umble prefer to order the
partitions chosen for the diagonal in a way such that any operations
defined on one side of the tensor product is reversed when performing
it on the other. However, for ease of programming, we shall adopt a
convention for internal representation in which the second argument is
stored in reverse. Nevertheless, our output functions take this
difference into consideration and make certain that everything
displayed for the user adheres to Saneblidze and Umble's convention.

A permutation corresponds to a pair of partitions by reading the first
partition as the falling sequences of the permutation, and the second
as the rising sequences. This corresponds to reading columns and rows
in the matrix representation used by Saneblidze and Umble. The matrix
constructed directly from a permutation is called a \emph{step
  matrix}, and we shall use this term to refer to the corresponding
partition pairs throughout. The code that follows below incorporates
this construction method.

Over fields of characteristic different from $2$, the terms of the
diagonal each have a sign. The code here developed will calculate the
sign along with the terms, and optionally strip the signs when
displaying the calculated faces. By the laziness properties of
Haskell, these calculations will only be performed as and when the
signs are requested.

\begingroup\par\noindent\advance\leftskip\mathindent\(
\begin{pboxed}\SaveRestoreHook
\column{B}{@{}l@{}}
\column{E}{@{}l@{}}
\>[B]{}\mathbf{type}\;\Conid{Face}\mathrel{=}(\Conid{Partition},\Conid{Partition}){}\<[E]
\\
\>[B]{}\mathbf{type}\;\Conid{SignFace}\mathrel{=}(\Conid{Int},\Conid{Partition},\Conid{Partition}){}\<[E]
\\[\blanklineskip]
\>[B]{}\Varid{stripSign}\mathbin{::}\Conid{SignFace}\to \Conid{Face}{}\<[E]
\\
\>[B]{}\Varid{stripSign}\;(\Varid{s},\Varid{p},\Varid{q})\mathrel{=}(\Varid{p},\Varid{q}){}\<[E]
\\[\blanklineskip]
\>[B]{}\Varid{buildFace}\mathbin{::}\Conid{Sequence}\to \Conid{SignFace}{}\<[E]
\\
\>[B]{}\Varid{buildFace}\;\Varid{p}\mathrel{=}\Varid{signFace}\;(\Varid{map}\;\Varid{sort}\;(\Varid{falling}\;\Varid{p}),\Varid{reverse}\;(\Varid{rising}\;\Varid{p})){}\<[E]
\ColumnHook
\end{pboxed}
\)\par\noindent\endgroup\resethooks

For an example consider the face denoted by Saneblidze and Umble as
$12|34|5\otimes 2|14|35$. This is a step matrix resulting from the
permutation $21435$. In the internal representation is would be
represented as \ensuremath{([\mskip1.5mu [\mskip1.5mu \mathrm{1},\mathrm{2}\mskip1.5mu],[\mskip1.5mu \mathrm{3},\mathrm{4}\mskip1.5mu],[\mskip1.5mu \mathrm{5}\mskip1.5mu]\mskip1.5mu],[\mskip1.5mu [\mskip1.5mu \mathrm{3},\mathrm{5}\mskip1.5mu],[\mskip1.5mu \mathrm{1},\mathrm{4}\mskip1.5mu],[\mskip1.5mu \mathrm{2}\mskip1.5mu]\mskip1.5mu])}. Without the
function call to \ensuremath{\Varid{sort}} in the first component, it would have been
\ensuremath{([\mskip1.5mu [\mskip1.5mu \mathrm{2},\mathrm{1}\mskip1.5mu],[\mskip1.5mu \mathrm{4},\mathrm{3}\mskip1.5mu],[\mskip1.5mu \mathrm{5}\mskip1.5mu]\mskip1.5mu],[\mskip1.5mu [\mskip1.5mu \mathrm{3},\mathrm{5}\mskip1.5mu],[\mskip1.5mu \mathrm{1},\mathrm{4}\mskip1.5mu],[\mskip1.5mu \mathrm{2}\mskip1.5mu]\mskip1.5mu])} and without the \ensuremath{\Varid{reverse}} in
the second component, it would have been
\ensuremath{([\mskip1.5mu [\mskip1.5mu \mathrm{1},\mathrm{2}\mskip1.5mu],[\mskip1.5mu \mathrm{3},\mathrm{4}\mskip1.5mu],[\mskip1.5mu \mathrm{5}\mskip1.5mu]\mskip1.5mu],[\mskip1.5mu [\mskip1.5mu \mathrm{2}\mskip1.5mu],[\mskip1.5mu \mathrm{3},\mathrm{5}\mskip1.5mu],[\mskip1.5mu \mathrm{1},\mathrm{4}\mskip1.5mu]\mskip1.5mu])}.

Thus, the benefit from the \ensuremath{\Varid{sort}} call is pure readability for the
user, whereas the benefit from the use of \ensuremath{\Varid{reverse}} is in ease of code
reuse, as will be seen when we implement the functions to generate
derived matrices.

The formula for the sign is complicated. It is helpful to view it as
the product of three different signs. One of these is calculated using
the sign of a permutation.  Suppose $\sigma$ is a permutation in
$S_{n+1}$ giving rise to the face indexed by $\mu'\otimes\mu$
with $\mu=\mu_r|\dots|\mu_1$. Then the concatenation of
the parts of $\mu$ retrieves the list of images for the partition
$\sigma$. We define $\operatorname{psgn}(\sigma)$ as the usual
permutation sign, and implement this as a simple orbit finder. The
code given here adjusts for the fact that lists in Haskell have
$0$-based indexing.

\begingroup\par\noindent\advance\leftskip\mathindent\(
\begin{pboxed}\SaveRestoreHook
\column{B}{@{}l@{}}
\column{3}{@{}l@{}}
\column{5}{@{}l@{}}
\column{7}{@{}l@{}}
\column{9}{@{}l@{}}
\column{11}{@{}l@{}}
\column{13}{@{}l@{}}
\column{14}{@{}l@{}}
\column{28}{@{}l@{}}
\column{E}{@{}l@{}}
\>[B]{}\Varid{orbit}\mathbin{::}\Conid{Int}\to \Conid{Sequence}\to \Conid{Sequence}{}\<[E]
\\
\>[B]{}\Varid{orbit}\;\Varid{a}\;\Varid{pi}\mathrel{=}\Varid{findOrbit}\;\Varid{a}\;[\mskip1.5mu \mskip1.5mu]{}\<[E]
\\
\>[B]{}\hsindent{3}{}\<[3]
\>[3]{}\mathbf{where}{}\<[E]
\\
\>[3]{}\hsindent{2}{}\<[5]
\>[5]{}\Varid{findOrbit}\;\Varid{a}\;\Varid{as}\mathrel{=}{}\<[E]
\\
\>[5]{}\hsindent{2}{}\<[7]
\>[7]{}\mathbf{if}\;{}\<[13]
\>[13]{}\Varid{a'}\in \Varid{as}{}\<[E]
\\
\>[5]{}\hsindent{2}{}\<[7]
\>[7]{}\mathbf{then}\;{}\<[13]
\>[13]{}(\Varid{sort}\mathbin{\circ}\Varid{nub})\;(\Varid{a}\mathbin{:}\Varid{as}){}\<[E]
\\
\>[5]{}\hsindent{2}{}\<[7]
\>[7]{}\mathbf{else}\;{}\<[13]
\>[13]{}\Varid{findOrbit}\;\Varid{a'}\;(\Varid{a}\mathbin{:}\Varid{as}){}\<[E]
\\
\>[7]{}\hsindent{4}{}\<[11]
\>[11]{}\mathbf{where}{}\<[E]
\\
\>[11]{}\hsindent{3}{}\<[14]
\>[14]{}\Varid{a'}\mathrel{=}\Varid{pi}\mathbin{!!}(\Varid{a}\mathbin{-}\mathrm{1}){}\<[E]
\\[\blanklineskip]
\>[B]{}\Varid{pSign}\mathbin{::}\Conid{Sequence}\to \Conid{Int}{}\<[E]
\\
\>[B]{}\Varid{pSign}\;\Varid{pi}\mathrel{=}\Varid{signPi}{}\<[E]
\\
\>[B]{}\hsindent{3}{}\<[3]
\>[3]{}\mathbf{where}{}\<[E]
\\
\>[3]{}\hsindent{2}{}\<[5]
\>[5]{}\Varid{getOrbits}\;\Varid{orbs}\;[\mskip1.5mu \mskip1.5mu]{}\<[28]
\>[28]{}\mathrel{=}\Varid{orbs}{}\<[E]
\\
\>[3]{}\hsindent{2}{}\<[5]
\>[5]{}\Varid{getOrbits}\;\Varid{orbs}\;(\Varid{p}\mathbin{:}\Varid{ps}){}\<[28]
\>[28]{}\mathrel{=}\Varid{getOrbits}\;(\Varid{o}\mathbin{:}\Varid{orbs})\;(\Varid{ps}\mathbin{\char92 \char92 }\Varid{o}){}\<[E]
\\
\>[5]{}\hsindent{2}{}\<[7]
\>[7]{}\mathbf{where}{}\<[E]
\\
\>[7]{}\hsindent{2}{}\<[9]
\>[9]{}\Varid{o}\mathrel{=}\Varid{orbit}\;\Varid{p}\;\Varid{pi}{}\<[E]
\\
\>[3]{}\hsindent{2}{}\<[5]
\>[5]{}\Varid{orbits}{}\<[28]
\>[28]{}\mathrel{=}\Varid{getOrbits}\;[\mskip1.5mu \mskip1.5mu]\;\Varid{pi}{}\<[E]
\\
\>[3]{}\hsindent{2}{}\<[5]
\>[5]{}\Varid{orbitLengths}{}\<[28]
\>[28]{}\mathrel{=}\Varid{map}\;\Varid{length}\;\Varid{orbits}{}\<[E]
\\
\>[3]{}\hsindent{2}{}\<[5]
\>[5]{}\Varid{evenCycles}{}\<[28]
\>[28]{}\mathrel{=}\Varid{filter}\;\Varid{even}\;\Varid{orbitLengths}{}\<[E]
\\
\>[3]{}\hsindent{2}{}\<[5]
\>[5]{}\Varid{signPi}{}\<[28]
\>[28]{}\mathrel{=}(\mathbin{-}\mathrm{1})\mathbin{\uparrow}(\Varid{length}\;\Varid{evenCycles}){}\<[E]
\ColumnHook
\end{pboxed}
\)\par\noindent\endgroup\resethooks

Furthermore, we let $\epsilon =
\sum_{i=1}^{r-1}i\cdot|\mu_i|$ and define the right-most
partition sign
$\operatorname{sgn}_r(\mu)=(-1)^\epsilon\operatorname{psgn}(\sigma)$.

\begingroup\par\noindent\advance\leftskip\mathindent\(
\begin{pboxed}\SaveRestoreHook
\column{B}{@{}l@{}}
\column{3}{@{}l@{}}
\column{5}{@{}l@{}}
\column{15}{@{}l@{}}
\column{E}{@{}l@{}}
\>[B]{}\Varid{signR}\mathbin{::}\Conid{Partition}\to \Conid{Int}{}\<[E]
\\
\>[B]{}\Varid{signR}\;\Varid{q}\mathrel{=}(\mathbin{-}\mathrm{1})\mathbin{\uparrow}\Varid{epsilon}\mathbin{*}(\Varid{pSign}\;\Varid{pi}){}\<[E]
\\
\>[B]{}\hsindent{3}{}\<[3]
\>[3]{}\mathbf{where}{}\<[E]
\\
\>[3]{}\hsindent{2}{}\<[5]
\>[5]{}\Varid{pi}{}\<[15]
\>[15]{}\mathrel{=}\Varid{concat}\;\Varid{q}{}\<[E]
\\
\>[3]{}\hsindent{2}{}\<[5]
\>[5]{}\Varid{epsilon}{}\<[15]
\>[15]{}\mathrel{=}\Varid{sum}\;\Varid{summands}{}\<[E]
\\
\>[3]{}\hsindent{2}{}\<[5]
\>[5]{}\Varid{summands}{}\<[15]
\>[15]{}\mathrel{=}\Varid{map}\;(\lambda \Varid{i}\to \Varid{i}\mathbin{*}(\Varid{qLengths}\mathbin{!!}(\Varid{i}\mathbin{-}\mathrm{1})))\;[\mskip1.5mu \mathrm{1}\mathinner{\ldotp\ldotp}((\Varid{length}\;\Varid{q})\mathbin{-}\mathrm{1})\mskip1.5mu]{}\<[E]
\\
\>[3]{}\hsindent{2}{}\<[5]
\>[5]{}\Varid{qLengths}{}\<[15]
\>[15]{}\mathrel{=}\Varid{map}\;\Varid{length}\;\Varid{q}{}\<[E]
\ColumnHook
\end{pboxed}
\)\par\noindent\endgroup\resethooks
  
Finally, we define the order-reversing permutation sign by
\[
\operatorname{orsgn}(\mu)=(-1)^{\frac12\left(\sum|U_i|^2 - (n+1)\right)}
\]

\begingroup\par\noindent\advance\leftskip\mathindent\(
\begin{pboxed}\SaveRestoreHook
\column{B}{@{}l@{}}
\column{3}{@{}l@{}}
\column{5}{@{}l@{}}
\column{20}{@{}l@{}}
\column{E}{@{}l@{}}
\>[B]{}\Varid{orSign}\mathbin{::}\Conid{Partition}\to \Conid{Int}{}\<[E]
\\
\>[B]{}\Varid{orSign}\;\Varid{p}\mathrel{=}(\mathbin{-}\mathrm{1})\mathbin{\uparrow}\Varid{exponent}{}\<[E]
\\
\>[B]{}\hsindent{3}{}\<[3]
\>[3]{}\mathbf{where}{}\<[E]
\\
\>[3]{}\hsindent{2}{}\<[5]
\>[5]{}\Varid{exponent}{}\<[20]
\>[20]{}\mathrel{=}\Varid{exponent2}\mathbin{\Varid{`div`}}\mathrm{2}{}\<[E]
\\
\>[3]{}\hsindent{2}{}\<[5]
\>[5]{}\Varid{exponent2}{}\<[20]
\>[20]{}\mathrel{=}(\Varid{sum}\;\Varid{lengthSquares})\mathbin{-}((\Varid{length}\mathbin{\circ}\Varid{concat})\;\Varid{p}){}\<[E]
\\
\>[3]{}\hsindent{2}{}\<[5]
\>[5]{}\Varid{lengthSquares}{}\<[20]
\>[20]{}\mathrel{=}\Varid{map}\;((\mathbin{\uparrow}\mathrm{2})\mathbin{\circ}\Varid{length})\;\Varid{p}{}\<[E]
\ColumnHook
\end{pboxed}
\)\par\noindent\endgroup\resethooks
  
Following Saneblidze-Umble \cite{saneblidze-umble-2004}, we now
allocate the sign 
\[
\operatorname{csgn}(p\otimes q)=(-1)^{|q|\choose2}
\operatorname{sgn}_r(p)\operatorname{orsgn}(q)
\]
 to a face
$p\otimes q$ given by a step matrix. We notice that in order to
reconcile the internal representation used in this program with the
 conventions in \cite{saneblidze-umble-2004}, we further need to
 reverse the order of the parts in $q$ to get the right sign.

\begingroup\par\noindent\advance\leftskip\mathindent\(
\begin{pboxed}\SaveRestoreHook
\column{B}{@{}l@{}}
\column{3}{@{}l@{}}
\column{5}{@{}l@{}}
\column{16}{@{}l@{}}
\column{26}{@{}l@{}}
\column{E}{@{}l@{}}
\>[B]{}\Varid{signFace}\mathbin{::}\Conid{Face}\to \Conid{SignFace}{}\<[E]
\\
\>[B]{}\Varid{signFace}\;(\Varid{p},\Varid{q})\mathrel{=}(\Varid{qSign}\mathbin{*}\Varid{rSign}\mathbin{*}\Varid{sign1},\Varid{p},\Varid{q}){}\<[E]
\\
\>[B]{}\hsindent{3}{}\<[3]
\>[3]{}\mathbf{where}{}\<[E]
\\
\>[3]{}\hsindent{2}{}\<[5]
\>[5]{}\Varid{qSign}{}\<[16]
\>[16]{}\mathrel{=}(\mathbin{-}\mathrm{1})\mathbin{\uparrow}\Varid{qExp}{}\<[E]
\\
\>[3]{}\hsindent{2}{}\<[5]
\>[5]{}\Varid{qExp}{}\<[16]
\>[16]{}\mathrel{=}(\Varid{choose2}\mathbin{\circ}\Varid{length})\;\Varid{q}{}\<[E]
\\
\>[3]{}\hsindent{2}{}\<[5]
\>[5]{}\Varid{rSign}{}\<[16]
\>[16]{}\mathrel{=}\Varid{orSign}\;{}\<[26]
\>[26]{}\Varid{p}{}\<[E]
\\
\>[3]{}\hsindent{2}{}\<[5]
\>[5]{}\Varid{sign1}{}\<[16]
\>[16]{}\mathrel{=}\Varid{signR}\;{}\<[26]
\>[26]{}(\Varid{reverse}\;\Varid{q}){}\<[E]
\\
\>[3]{}\hsindent{2}{}\<[5]
\>[5]{}\Varid{choose2}\;\Varid{n}{}\<[16]
\>[16]{}\mathrel{=}\Varid{n}\mathbin{*}(\Varid{n}\mathbin{-}\mathrm{1})\mathbin{\Varid{`div`}}\mathrm{2}{}\<[E]
\ColumnHook
\end{pboxed}
\)\par\noindent\endgroup\resethooks

Once the calculations have been performed, the user is likely to wish
for a readable display of the calculated faces. In order to do this,
the following code gives several functions that display faces of the
permutahedron in a way familiar to the reader.

\begingroup\par\noindent\advance\leftskip\mathindent\(
\begin{pboxed}\SaveRestoreHook
\column{B}{@{}l@{}}
\column{3}{@{}l@{}}
\column{5}{@{}l@{}}
\column{7}{@{}l@{}}
\column{11}{@{}l@{}}
\column{13}{@{}l@{}}
\column{16}{@{}l@{}}
\column{19}{@{}l@{}}
\column{E}{@{}l@{}}
\>[B]{}\Varid{showSignFace}\mathbin{::}\Conid{SignFace}\to \Conid{String}{}\<[E]
\\
\>[B]{}\Varid{showSignFace}\;\Varid{f}\mathord{@}(\Varid{s},\anonymous ,\anonymous )\mathrel{=}{}\<[E]
\\
\>[B]{}\hsindent{3}{}\<[3]
\>[3]{}\mathbf{case}\;\Varid{s}\;\mathbf{of}{}\<[E]
\\
\>[B]{}\hsindent{3}{}\<[3]
\>[3]{}\mathrm{1}{}\<[7]
\>[7]{}\to \text{\tt \char34 +\char34}\plus (\Varid{showFace}\mathbin{\circ}\Varid{stripSign})\;\Varid{f}{}\<[E]
\\
\>[B]{}\hsindent{3}{}\<[3]
\>[3]{}\mathbin{-}\mathrm{1}{}\<[7]
\>[7]{}\to \text{\tt \char34 -\char34}\plus (\Varid{showFace}\mathbin{\circ}\Varid{stripSign})\;\Varid{f}{}\<[E]
\\
\>[B]{}\hsindent{3}{}\<[3]
\>[3]{}\mathrm{0}{}\<[7]
\>[7]{}\to \text{\tt \char34 \char34}{}\<[E]
\\
\>[B]{}\hsindent{3}{}\<[3]
\>[3]{}\Varid{a}{}\<[7]
\>[7]{}\to (\Varid{show}\;\Varid{a})\plus \text{\tt \char34 .\char34}\plus (\Varid{showFace}\mathbin{\circ}\Varid{stripSign})\;\Varid{f}{}\<[E]
\\[\blanklineskip]
\>[B]{}\Varid{showFace}{}\<[11]
\>[11]{}\mathbin{::}\Conid{Face}\to \Conid{String}{}\<[E]
\\
\>[B]{}\Varid{showFace}\mathrel{=}\Varid{showFaceTemplate}\;\Varid{showPartition}{}\<[E]
\\[\blanklineskip]
\>[B]{}\Varid{showFaceShort}\mathbin{::}\Conid{Face}\to \Conid{String}{}\<[E]
\\
\>[B]{}\Varid{showFaceShort}\mathrel{=}\Varid{showFaceTemplate}\;\Varid{showPartitionShort}{}\<[E]
\\[\blanklineskip]
\>[B]{}\Varid{showFaceTemplate}\mathbin{::}(\Conid{Partition}\to \Conid{String})\to \Conid{Face}\to \Conid{String}{}\<[E]
\\
\>[B]{}\Varid{showFaceTemplate}\;\Varid{showP}\;(\Varid{u},\Varid{v})\mathrel{=}\Varid{showP}\;\Varid{u}\plus \text{\tt \char34 x\char34}\plus (\Varid{showP}\mathbin{\circ}\Varid{reverse})\;\Varid{v}{}\<[E]
\\[\blanklineskip]
\>[B]{}\Varid{showPartitionShort}\mathbin{::}\Conid{Partition}\to \Conid{String}{}\<[E]
\\
\>[B]{}\Varid{showPartitionShort}\mathrel{=}\Varid{filter}\;(\not\equiv \text{\tt ','})\mathbin{\circ}\Varid{showPartition}{}\<[E]
\\[\blanklineskip]
\>[B]{}\Varid{showPartition}\mathbin{::}\Conid{Partition}\to \Conid{String}{}\<[E]
\\
\>[B]{}\Varid{showPartition}\;\Varid{p}\mathrel{=}\Varid{pString}\;\Varid{p}{}\<[E]
\\
\>[B]{}\hsindent{3}{}\<[3]
\>[3]{}\mathbf{where}{}\<[E]
\\
\>[3]{}\hsindent{2}{}\<[5]
\>[5]{}\Varid{pString}{}\<[19]
\>[19]{}\mathrel{=}\Varid{concat}\mathbin{\circ}\Varid{intersperse}\;\text{\tt \char34 |\char34}\mathbin{\circ}\Varid{partsStrings}{}\<[E]
\\
\>[3]{}\hsindent{2}{}\<[5]
\>[5]{}\Varid{partsStrings}{}\<[19]
\>[19]{}\mathrel{=}\Varid{map}\;(\Varid{concat}\mathbin{\circ}\Varid{intersperse}\;\text{\tt \char34 ,\char34}\mathbin{\circ}\Varid{map}\;\Varid{show}){}\<[E]
\\[\blanklineskip]
\>[B]{}\Varid{showMatrix}\mathbin{::}\Conid{Face}\to \Conid{String}{}\<[E]
\\
\>[B]{}\Varid{showMatrix}\;(\Varid{f1},\Varid{f2})\mathrel{=}\Varid{unlines}\mathbin{\$}\Varid{map}\;(\Varid{showLine}\;\Varid{f1})\;\Varid{f2}{}\<[E]
\\
\>[B]{}\hsindent{3}{}\<[3]
\>[3]{}\mathbf{where}{}\<[E]
\\
\>[3]{}\hsindent{2}{}\<[5]
\>[5]{}\Varid{showLine}\;{}\<[16]
\>[16]{}\Varid{a}\;\Varid{b}\mathrel{=}\Varid{concatMap}\;(\Varid{flip}\;\Varid{showPoint}\;\Varid{b})\;\Varid{a}{}\<[E]
\\
\>[3]{}\hsindent{2}{}\<[5]
\>[5]{}\Varid{showPoint}\;{}\<[16]
\>[16]{}\Varid{a}\;\Varid{b}\mathrel{=}{}\<[E]
\\
\>[5]{}\hsindent{2}{}\<[7]
\>[7]{}\mathbf{if}\;{}\<[13]
\>[13]{}\Varid{intersect}\;\Varid{a}\;\Varid{b}\not\equiv [\mskip1.5mu \mskip1.5mu]{}\<[E]
\\
\>[5]{}\hsindent{2}{}\<[7]
\>[7]{}\mathbf{then}\;{}\<[13]
\>[13]{}\Varid{show}\mathbin{\$}\Varid{head}\mathbin{\$}\Varid{intersect}\;\Varid{a}\;\Varid{b}{}\<[E]
\\
\>[5]{}\hsindent{2}{}\<[7]
\>[7]{}\mathbf{else}\;{}\<[13]
\>[13]{}\text{\tt \char34 .\char34}{}\<[E]
\ColumnHook
\end{pboxed}
\)\par\noindent\endgroup\resethooks
%$

Using this code, we can gain string representations of the face
$+12|34|5\otimes 2|14|35$ as follows.

\begin{tabular}{cc|cc}
  \ensuremath{\Varid{showSignFace}}  & \text{\tt \char43{}1\char44{}2\char124{}3\char44{}4\char124{}5x2\char124{}1\char44{}4\char124{}3\char44{}5} & & \text{\tt \char46{}35} \\
  \ensuremath{\Varid{showFace}}      & \text{\tt 1\char44{}2\char124{}3\char44{}4\char124{}5x2\char124{}1\char44{}4\char124{}3\char44{}5}  & \ensuremath{\Varid{showMatrix}} & \text{\tt 14\char46{}} \\
  \ensuremath{\Varid{showFaceShort}} & \text{\tt 12\char124{}34\char124{}5x2\char124{}14\char124{}35}  & & \text{\tt 2\char46{}\char46{}} \\
\end{tabular}

We shall need to enumerate all step matrices, and thus, we need a
function to enumerate all permutations in $S_n$. The following code
takes care of this.

\begingroup\par\noindent\advance\leftskip\mathindent\(
\begin{pboxed}\SaveRestoreHook
\column{B}{@{}l@{}}
\column{3}{@{}l@{}}
\column{5}{@{}l@{}}
\column{22}{@{}l@{}}
\column{E}{@{}l@{}}
\>[B]{}\Varid{permutations}\mathbin{::}\Conid{Int}\to [\mskip1.5mu \Conid{Sequence}\mskip1.5mu]{}\<[E]
\\
\>[B]{}\Varid{permutations}\;\Varid{n}\mathrel{=}\Varid{permuteList}\;[\mskip1.5mu \mathrm{1}\mathinner{\ldotp\ldotp}\Varid{n}\mskip1.5mu]{}\<[E]
\\
\>[B]{}\hsindent{3}{}\<[3]
\>[3]{}\mathbf{where}{}\<[E]
\\
\>[3]{}\hsindent{2}{}\<[5]
\>[5]{}\Varid{permuteList}\;[\mskip1.5mu \mskip1.5mu]{}\<[22]
\>[22]{}\mathrel{=}[\mskip1.5mu \mskip1.5mu]{}\<[E]
\\
\>[3]{}\hsindent{2}{}\<[5]
\>[5]{}\Varid{permuteList}\;[\mskip1.5mu \Varid{a}\mskip1.5mu]{}\<[22]
\>[22]{}\mathrel{=}[\mskip1.5mu [\mskip1.5mu \Varid{a}\mskip1.5mu]\mskip1.5mu]{}\<[E]
\\
\>[3]{}\hsindent{2}{}\<[5]
\>[5]{}\Varid{permuteList}\;\Varid{l}{}\<[22]
\>[22]{}\mathrel{=}\Varid{concatMap}\;(\lambda \Varid{x}\to \Varid{map}\;(\Varid{x}\mathbin{:})\;(\Varid{permuteList}\;(\Varid{l}\mathbin{\char92 \char92 }[\mskip1.5mu \Varid{x}\mskip1.5mu])))\;\Varid{l}{}\<[E]
\ColumnHook
\end{pboxed}
\)\par\noindent\endgroup\resethooks

From these matrices, we then generate the closure under two
operations: first downshift and then rightshift. These are, with the
representation I have chosen, the same operation on either side of the
tensor product. So it is sufficient to define a single operation on a
partition.

\begin{definition} 
  Let $\pi=\pi_1|\pi_2|\dots|\pi_k$ be a partition.  A proper
  subset $M$ of some part $\pi_j$ in $\pi$ with $j<k$ is
  \emph{admissible} with respect to a partition
  $\mu=\mu_1|\mu_2|\dots|\mu_r$ if $\min
  M>\max\pi_{j+1}$ and, supposing that $\min M$ occurs in $\mu_k$,
  then $\pi_{j+1}\cap\bigcup_{t=k}^r \mu_r=\emptyset$.

  An $M$-shift of a partition $\pi$ with respect to a partition
  $\mu$, with $M$ an admissible subset of $\pi_j$ with respect to
  $\mu$, is the operation that returns the new partition
  $\pi_1|\dots|\pi_j\setminus M|\pi_{j+1}\cup M|\dots|\pi_k$.
\end{definition}

Note that we only need to define the shift operations on the first
factor of a face, since we can always conjugate it with the twist
operation that interchanges factors in a tensor product.

Admissibility is not enough for a particular move to be
permitted. There is also a condition on the sequence of moves that led
to the partition $\pi$: we require that the moves
progress through the matrix, so that if we have performed a move on an
admissible subset $M\subset\pi_i$, then all subsequent moves have
to occur on subsets $M'\subset\pi_j$ with $j>i$. 

The conditions placed on permissibility apart from the admissibility
condition ensure that no matrix occurs as a derived matrix from more
than one step matrix. This makes it clear
that the process of generating new derived matrices will stop at some
point. Thus, among the moves \emph{not} permitted we find
\[
\begin{pmatrix}
  & 1 & 2 \\
  4 & 5 & \\
  3 & & 
\end{pmatrix}
\Rightarrow
\begin{pmatrix}
  & 1 & 2 \\
  4 & & 5 \\
  3 & & 
\end{pmatrix}
\Rightarrow
\begin{pmatrix}
  & 1 & 2 \\
  & 4 & 5 \\
  3 & & 
\end{pmatrix}
\]

Checking admissibility for a certain subset of $\pi_1$ is just a
matter of constructing programmatic checks for all conditions in the
definition, and then ensuring that all these conditions hold.

Since Haskell evaluates as lazily as possible we can get around the
fact that certain conditions can only be defined once earlier
conditions are seen to hold, since the resulting code will break
evaluation as soon as a non-holding condition has been found. Thus,
the latter conditions are not evaluated unless the earlier have been
found to hold, thus making the latter conditions well-defined. 

\begingroup\par\noindent\advance\leftskip\mathindent\(
\begin{pboxed}\SaveRestoreHook
\column{B}{@{}l@{}}
\column{3}{@{}l@{}}
\column{5}{@{}l@{}}
\column{20}{@{}l@{}}
\column{E}{@{}l@{}}
\>[B]{}\Varid{isAdmissible}\mathbin{::}\Conid{SignFace}\to [\mskip1.5mu \Conid{Int}\mskip1.5mu]\to \Conid{Bool}{}\<[E]
\\
\>[B]{}\Varid{isAdmissible}\;\Varid{f}\mathord{@}(\anonymous ,\Varid{pi},\Varid{mu})\;\Varid{m}\mathrel{=}\Varid{admitted}{}\<[E]
\\
\>[B]{}\hsindent{3}{}\<[3]
\>[3]{}\mathbf{where}{}\<[E]
\\
\>[3]{}\hsindent{2}{}\<[5]
\>[5]{}\Varid{mIntersectPi}{}\<[20]
\>[20]{}\mathrel{=}\Varid{map}\;(\Varid{intersect}\;\Varid{m})\;\Varid{pi}{}\<[E]
\\
\>[3]{}\hsindent{2}{}\<[5]
\>[5]{}\Varid{partsWithM}{}\<[20]
\>[20]{}\mathrel{=}\Varid{findIndices}\;((\equiv \Varid{length}\;\Varid{m})\mathbin{\circ}\Varid{length})\;\Varid{mIntersectPi}{}\<[E]
\\
\>[3]{}\hsindent{2}{}\<[5]
\>[5]{}\Varid{mInUniquePart}{}\<[20]
\>[20]{}\mathrel{=}(\mathrm{1}\equiv \Varid{length}\;\Varid{partsWithM}){}\<[E]
\\
\>[3]{}\hsindent{2}{}\<[5]
\>[5]{}\Varid{j}{}\<[20]
\>[20]{}\mathrel{=}\Varid{head}\;\Varid{partsWithM}{}\<[E]
\\
\>[3]{}\hsindent{2}{}\<[5]
\>[5]{}\Varid{jLessK}{}\<[20]
\>[20]{}\mathrel{=}\Varid{j}\mathbin{<}\Varid{length}\;\Varid{pi}{}\<[E]
\\
\>[3]{}\hsindent{2}{}\<[5]
\>[5]{}\Varid{pi\char95 j}{}\<[20]
\>[20]{}\mathrel{=}\Varid{pi}\mathbin{!!}\Varid{j}{}\<[E]
\\
\>[3]{}\hsindent{2}{}\<[5]
\>[5]{}\Varid{pi\char95 j1}{}\<[20]
\>[20]{}\mathrel{=}\Varid{pi}\mathbin{!!}(\Varid{j}\mathbin{+}\mathrm{1}){}\<[E]
\\
\>[3]{}\hsindent{2}{}\<[5]
\>[5]{}\Varid{properSubset}{}\<[20]
\>[20]{}\mathrel{=}\Varid{mInUniquePart}\mathrel{\wedge}\Varid{jLessK}\mathrel{\wedge}(\Varid{length}\;\Varid{m}\mathbin{<}\Varid{length}\;\Varid{pi\char95 j}){}\<[E]
\\
\>[3]{}\hsindent{2}{}\<[5]
\>[5]{}\Varid{minLargerMax}{}\<[20]
\>[20]{}\mathrel{=}(\Varid{minimum}\;\Varid{m}\mathbin{>}\Varid{maximum}\;\Varid{pi\char95 j1}){}\<[E]
\\
\>[3]{}\hsindent{2}{}\<[5]
\>[5]{}\Varid{k}{}\<[20]
\>[20]{}\mathrel{=}\Varid{fromJust}\;(\Varid{findIndex}\;(\Varid{minimum}\;\Varid{m}\in )\;\Varid{mu}){}\<[E]
\\
\>[3]{}\hsindent{2}{}\<[5]
\>[5]{}\Varid{mus}{}\<[20]
\>[20]{}\mathrel{=}\Varid{concat}\;(\Varid{drop}\;\Varid{k}\;\Varid{mu}){}\<[E]
\\
\>[3]{}\hsindent{2}{}\<[5]
\>[5]{}\Varid{allzero}{}\<[20]
\>[20]{}\mathrel{=}\Varid{null}\;(\Varid{intersect}\;\Varid{pi\char95 j1}\;\Varid{mus}){}\<[E]
\\
\>[3]{}\hsindent{2}{}\<[5]
\>[5]{}\Varid{admitted}{}\<[20]
\>[20]{}\mathrel{=}\Varid{properSubset}\mathrel{\wedge}\Varid{minLargerMax}\mathrel{\wedge}\Varid{allzero}{}\<[E]
\ColumnHook
\end{pboxed}
\)\par\noindent\endgroup\resethooks

At this point, moving an admissible subset is a matter of assembling a
new list with the appropriate parts modified to remove the set from
one part and adjoin it to the next, and adjusting the sign
accordingly. In order to figure out the sign, we thus
need to move our elements in increasing order, one by one, and catch
the sign changes as we do.

Suppose, following Saneblidze and Umble \cite{saneblidze-umble-2004},
that we move the single element $x\in\pi_i$ to $\pi_{i+1}$,
and our face is currently represented as $\pi\otimes\mu$. We
define the upper and lower cuts by $(a,S]=\{s\in S|s>a\}$ and
$[S,a)=\{s\in S|s<a\}$ respectively. Then
\[
\operatorname{csgn}(R_x\pi\otimes\mu)=
-\operatorname{csgn}(\pi\otimes\mu)\cdot
(-1)^{|(x,\pi_i]\cup[\pi_{i+1},x)|}
\]
Down moves are simply right moves conjugated with a transposition, or
in other words, this operation conjugated with the twist map $a\otimes
b\mapsto b\otimes a$. Thus, we need only handle this case, since the
other case follows symmetrically.

The function here written handles non-admissible subsets gracefully.

\begingroup\par\noindent\advance\leftskip\mathindent\(
\begin{pboxed}\SaveRestoreHook
\column{B}{@{}l@{}}
\column{3}{@{}l@{}}
\column{5}{@{}l@{}}
\column{7}{@{}l@{}}
\column{9}{@{}l@{}}
\column{11}{@{}l@{}}
\column{13}{@{}l@{}}
\column{19}{@{}l@{}}
\column{21}{@{}l@{}}
\column{E}{@{}l@{}}
\>[B]{}\Varid{moveSubset}\mathbin{::}\Conid{SignFace}\to [\mskip1.5mu \Conid{Int}\mskip1.5mu]\to \Conid{Maybe}\;\Conid{SignFace}{}\<[E]
\\
\>[B]{}\Varid{moveSubset}\;\Varid{f}\mathord{@}(\anonymous ,\Varid{p'},\anonymous )\;\Varid{m}\mathrel{=}{}\<[E]
\\
\>[B]{}\hsindent{3}{}\<[3]
\>[3]{}\mathbf{if}\;{}\<[9]
\>[9]{}\Varid{isAdmissible}\;\Varid{f}\;\Varid{m}{}\<[E]
\\
\>[B]{}\hsindent{3}{}\<[3]
\>[3]{}\mathbf{then}\;{}\<[9]
\>[9]{}\Conid{Just}\;(\Varid{foldl'}\;\Varid{moveElement}\;\Varid{f}\;(\Varid{sort}\;\Varid{m})){}\<[E]
\\
\>[B]{}\hsindent{3}{}\<[3]
\>[3]{}\mathbf{else}\;{}\<[9]
\>[9]{}\Conid{Nothing}{}\<[E]
\\
\>[3]{}\hsindent{2}{}\<[5]
\>[5]{}\mathbf{where}{}\<[E]
\\
\>[5]{}\hsindent{2}{}\<[7]
\>[7]{}\Varid{moveElement}\mathbin{::}\Conid{SignFace}\to \Conid{Int}\to \Conid{SignFace}{}\<[E]
\\
\>[5]{}\hsindent{2}{}\<[7]
\>[7]{}\Varid{moveElement}\;(\Varid{s},\Varid{p},\Varid{q})\;\Varid{e}\mathrel{=}(\Varid{s'},\Varid{p'},\Varid{q}){}\<[E]
\\
\>[7]{}\hsindent{2}{}\<[9]
\>[9]{}\mathbf{where}{}\<[E]
\\
\>[9]{}\hsindent{2}{}\<[11]
\>[11]{}(\Conid{Just}\;\Varid{i}){}\<[21]
\>[21]{}\mathrel{=}\Varid{findIndex}\;(\Varid{e}\in )\;\Varid{p}{}\<[E]
\\
\>[9]{}\hsindent{2}{}\<[11]
\>[11]{}\Varid{pi}{}\<[21]
\>[21]{}\mathrel{=}\Varid{p}\mathbin{!!}\Varid{i}{}\<[E]
\\
\>[9]{}\hsindent{2}{}\<[11]
\>[11]{}\Varid{pi1}{}\<[21]
\>[21]{}\mathrel{=}\Varid{p}\mathbin{!!}(\Varid{i}\mathbin{+}\mathrm{1}){}\<[E]
\\
\>[9]{}\hsindent{2}{}\<[11]
\>[11]{}\Varid{pmoved}{}\<[21]
\>[21]{}\mathrel{=}(\Varid{take}\;\Varid{i}\;\Varid{p})\plus [\mskip1.5mu \Varid{pi}\mathbin{\char92 \char92 }[\mskip1.5mu \Varid{e}\mskip1.5mu],\Varid{pi1}\plus [\mskip1.5mu \Varid{e}\mskip1.5mu]\mskip1.5mu]\plus (\Varid{drop}\;(\Varid{i}\mathbin{+}\mathrm{2})\;\Varid{p}){}\<[E]
\\
\>[9]{}\hsindent{2}{}\<[11]
\>[11]{}\Varid{lowercut}{}\<[21]
\>[21]{}\mathrel{=}\Varid{filter}\;(\mathbin{>}\Varid{e})\;\Varid{pi}{}\<[E]
\\
\>[9]{}\hsindent{2}{}\<[11]
\>[11]{}\Varid{uppercut}{}\<[21]
\>[21]{}\mathrel{=}\Varid{filter}\;(\mathbin{<}\Varid{e})\;\Varid{pi1}{}\<[E]
\\
\>[9]{}\hsindent{2}{}\<[11]
\>[11]{}\Varid{expmoved}{}\<[21]
\>[21]{}\mathrel{=}\Varid{length}\;(\Varid{lowercut}\plus \Varid{uppercut}){}\<[E]
\\
\>[9]{}\hsindent{2}{}\<[11]
\>[11]{}(\Varid{s'},\Varid{p'}){}\<[21]
\>[21]{}\mathrel{=}{}\<[E]
\\
\>[11]{}\hsindent{2}{}\<[13]
\>[13]{}\mathbf{if}\;{}\<[19]
\>[19]{}\Varid{e}\mathbin{>}\Varid{maximum}\;\Varid{pi1}{}\<[E]
\\
\>[11]{}\hsindent{2}{}\<[13]
\>[13]{}\mathbf{then}\;{}\<[19]
\>[19]{}(\mathbin{-}\Varid{s}\mathbin{*}(\mathbin{-}\mathrm{1})\mathbin{\uparrow}\Varid{expmoved},\Varid{pmoved}){}\<[E]
\\
\>[11]{}\hsindent{2}{}\<[13]
\>[13]{}\mathbf{else}\;{}\<[19]
\>[19]{}(\Varid{s},\Varid{p}){}\<[E]
\ColumnHook
\end{pboxed}
\)\par\noindent\endgroup\resethooks
%$

Now, to construct the Saneblidze-Umble diagonal as the closure of the
step matrices under these shift operations, we begin by enumerating
all admissible subsets of a single partition. We need only to check
such sets for admissibility that are proper subsets of some part and
consist of elements larger than the maximal element in the next part,
since otherwise they would not fulfill even the first few conditions.
\begingroup\par\noindent\advance\leftskip\mathindent\(
\begin{pboxed}\SaveRestoreHook
\column{B}{@{}l@{}}
\column{3}{@{}l@{}}
\column{5}{@{}l@{}}
\column{7}{@{}l@{}}
\column{9}{@{}l@{}}
\column{11}{@{}l@{}}
\column{17}{@{}l@{}}
\column{20}{@{}l@{}}
\column{27}{@{}l@{}}
\column{32}{@{}l@{}}
\column{E}{@{}l@{}}
\>[B]{}\Varid{admissiblesInPin}\mathbin{::}\Conid{Int}\to \Conid{SignFace}\to [\mskip1.5mu [\mskip1.5mu \Conid{Int}\mskip1.5mu]\mskip1.5mu]{}\<[E]
\\
\>[B]{}\Varid{admissiblesInPin}\;\Varid{i}\;\Varid{f}\mathord{@}(\anonymous ,\Varid{pi},\anonymous ){}\<[32]
\>[32]{}\mathrel{=}\Varid{filter}\;(\neg \mathbin{\circ}\Varid{null})\;\Varid{admissibleSets}{}\<[E]
\\
\>[B]{}\hsindent{3}{}\<[3]
\>[3]{}\mathbf{where}{}\<[E]
\\
\>[3]{}\hsindent{2}{}\<[5]
\>[5]{}\Varid{admissibleSets}\mathrel{=}{}\<[E]
\\
\>[5]{}\hsindent{2}{}\<[7]
\>[7]{}\mathbf{if}\;\Varid{i}\mathbin{+}\mathrm{2}\mathbin{>}\Varid{length}\;\Varid{pi}{}\<[E]
\\
\>[5]{}\hsindent{2}{}\<[7]
\>[7]{}\mathbf{then}\;[\mskip1.5mu \mskip1.5mu]{}\<[E]
\\
\>[5]{}\hsindent{2}{}\<[7]
\>[7]{}\mathbf{else}\;\Varid{filter}\;{}\<[20]
\>[20]{}(\Varid{isAdmissible}\;\Varid{f})\;\Varid{candidates}{}\<[E]
\\
\>[7]{}\hsindent{2}{}\<[9]
\>[9]{}\mathbf{where}{}\<[E]
\\
\>[9]{}\hsindent{2}{}\<[11]
\>[11]{}\Varid{candidates}{}\<[27]
\>[27]{}\mathrel{=}\Varid{filter}\;(\neg \mathbin{\circ}\Varid{null})\;(\Varid{subsets}\;\Varid{large}){}\<[E]
\\
\>[9]{}\hsindent{2}{}\<[11]
\>[11]{}\Varid{large}{}\<[27]
\>[27]{}\mathrel{=}\Varid{filter}\;(\mathbin{>}\Varid{m})\;(\Varid{pi}\mathbin{!!}\Varid{i}){}\<[E]
\\
\>[9]{}\hsindent{2}{}\<[11]
\>[11]{}\Varid{m}{}\<[27]
\>[27]{}\mathrel{=}\Varid{maximum}\;(\Varid{pi}\mathbin{!!}(\Varid{i}\mathbin{+}\mathrm{1})){}\<[E]
\\[\blanklineskip]
\>[B]{}\Varid{subsets}\mathbin{::}[\mskip1.5mu \Varid{a}\mskip1.5mu]\to [\mskip1.5mu [\mskip1.5mu \Varid{a}\mskip1.5mu]\mskip1.5mu]{}\<[E]
\\
\>[B]{}\Varid{subsets}\;[\mskip1.5mu \mskip1.5mu]{}\<[17]
\>[17]{}\mathrel{=}[\mskip1.5mu [\mskip1.5mu \mskip1.5mu]\mskip1.5mu]{}\<[E]
\\
\>[B]{}\Varid{subsets}\;(\Varid{a}\mathbin{:}\Varid{as}){}\<[17]
\>[17]{}\mathrel{=}\Varid{map}\;(\Varid{a}\mathbin{:})\;(\Varid{subsets}\;\Varid{as})\plus \Varid{subsets}\;\Varid{as}{}\<[E]
\ColumnHook
\end{pboxed}
\)\par\noindent\endgroup\resethooks

Now that we can enumerate all admissible subsets in a partition, it is
a simple matter to take a face and form all derived faces by
generating more and more admissible subsets and moving these. We
handle the moves in the second partition, as mentioned before, by
conjugating the shift operation by the twist map.
Valid derived faces are such that are reached by first moving subsets
right, and then moving them down. 
\begingroup\par\noindent\advance\leftskip\mathindent\(
\begin{pboxed}\SaveRestoreHook
\column{B}{@{}l@{}}
\column{3}{@{}l@{}}
\column{5}{@{}l@{}}
\column{7}{@{}l@{}}
\column{9}{@{}l@{}}
\column{13}{@{}l@{}}
\column{23}{@{}l@{}}
\column{E}{@{}l@{}}
\>[B]{}\Varid{twist}\mathbin{::}\Conid{SignFace}\to \Conid{SignFace}{}\<[E]
\\
\>[B]{}\Varid{twist}\;(\Varid{s},\Varid{a},\Varid{b})\mathrel{=}(\Varid{s},\Varid{b},\Varid{a}){}\<[E]
\\[\blanklineskip]
\>[B]{}\Varid{derivedFaces}\mathbin{::}\Conid{SignFace}\to [\mskip1.5mu \Conid{SignFace}\mskip1.5mu]{}\<[E]
\\
\>[B]{}\Varid{derivedFaces}\;\Varid{f}\mathord{@}(\anonymous ,\Varid{p},\Varid{q})\mathrel{=}\Varid{derivedRightQ}\;\mathrm{0}\;[\mskip1.5mu \mskip1.5mu]\;[\mskip1.5mu \Varid{f}\mskip1.5mu]{}\<[E]
\\
\>[B]{}\hsindent{3}{}\<[3]
\>[3]{}\mathbf{where}{}\<[E]
\\
\>[3]{}\hsindent{2}{}\<[5]
\>[5]{}\Varid{lp}\mathrel{=}\Varid{length}\;\Varid{p}{}\<[E]
\\
\>[3]{}\hsindent{2}{}\<[5]
\>[5]{}\Varid{lq}\mathrel{=}\Varid{length}\;\Varid{q}{}\<[E]
\\
\>[3]{}\hsindent{2}{}\<[5]
\>[5]{}\Varid{derivedRightQ}\;\Varid{i}\;\Varid{r}\;[\mskip1.5mu \mskip1.5mu]\mathrel{=}{}\<[E]
\\
\>[5]{}\hsindent{2}{}\<[7]
\>[7]{}\mathbf{if}\;{}\<[13]
\>[13]{}\Varid{i}\geq \Varid{lp}\mathbin{-}\mathrm{2}{}\<[E]
\\
\>[5]{}\hsindent{2}{}\<[7]
\>[7]{}\mathbf{then}\;{}\<[13]
\>[13]{}\Varid{derivedDownQ}\;\mathrm{0}\;[\mskip1.5mu \mskip1.5mu]\;\Varid{r}{}\<[E]
\\
\>[5]{}\hsindent{2}{}\<[7]
\>[7]{}\mathbf{else}\;{}\<[13]
\>[13]{}\Varid{derivedRightQ}\;(\Varid{i}\mathbin{+}\mathrm{1})\;[\mskip1.5mu \mskip1.5mu]\;\Varid{r}{}\<[E]
\\
\>[3]{}\hsindent{2}{}\<[5]
\>[5]{}\Varid{derivedRightQ}\;\Varid{i}\;\Varid{r}\;(\Varid{s}\mathbin{:}\Varid{ss})\mathrel{=}\Varid{derivedRightQ}\;\Varid{i}\;(\Varid{s}\mathbin{:}\Varid{r})\;(\Varid{ss}\plus \Varid{rights}){}\<[E]
\\
\>[5]{}\hsindent{2}{}\<[7]
\>[7]{}\mathbf{where}{}\<[E]
\\
\>[7]{}\hsindent{2}{}\<[9]
\>[9]{}\Varid{rights}\mathrel{=}\Varid{mapMaybe}\;(\Varid{moveSubset}\;\Varid{s})\;(\Varid{admissiblesInPin}\;\Varid{i}\;\Varid{s}){}\<[E]
\\
\>[3]{}\hsindent{2}{}\<[5]
\>[5]{}\Varid{derivedDownQ}\;\Varid{i}\;\Varid{d}\;[\mskip1.5mu \mskip1.5mu]\mathrel{=}{}\<[E]
\\
\>[5]{}\hsindent{2}{}\<[7]
\>[7]{}\mathbf{if}\;{}\<[13]
\>[13]{}\Varid{i}\geq \Varid{lq}\mathbin{-}\mathrm{2}{}\<[E]
\\
\>[5]{}\hsindent{2}{}\<[7]
\>[7]{}\mathbf{then}\;{}\<[13]
\>[13]{}\Varid{d}{}\<[E]
\\
\>[5]{}\hsindent{2}{}\<[7]
\>[7]{}\mathbf{else}\;{}\<[13]
\>[13]{}\Varid{derivedDownQ}\;(\Varid{i}\mathbin{+}\mathrm{1})\;[\mskip1.5mu \mskip1.5mu]\;\Varid{d}{}\<[E]
\\
\>[3]{}\hsindent{2}{}\<[5]
\>[5]{}\Varid{derivedDownQ}\;\Varid{i}\;\Varid{d}\;(\Varid{s}\mathbin{:}\Varid{ss})\mathrel{=}\Varid{derivedDownQ}\;\Varid{i}\;(\Varid{s}\mathbin{:}\Varid{d})\;(\Varid{ss}\plus \Varid{downs}){}\<[E]
\\
\>[5]{}\hsindent{2}{}\<[7]
\>[7]{}\mathbf{where}{}\<[E]
\\
\>[7]{}\hsindent{2}{}\<[9]
\>[9]{}\Varid{s'}{}\<[23]
\>[23]{}\mathrel{=}\Varid{twist}\;\Varid{s}{}\<[E]
\\
\>[7]{}\hsindent{2}{}\<[9]
\>[9]{}\Varid{downs}{}\<[23]
\>[23]{}\mathrel{=}\Varid{map}\;\Varid{twist}\;\Varid{twistedDowns}{}\<[E]
\\
\>[7]{}\hsindent{2}{}\<[9]
\>[9]{}\Varid{twistedDowns}{}\<[23]
\>[23]{}\mathrel{=}(\Varid{mapMaybe}\;(\Varid{moveSubset}\;\Varid{s'})\;(\Varid{admissiblesInPin}\;\Varid{i}\;\Varid{s'})){}\<[E]
\ColumnHook
\end{pboxed}
\)\par\noindent\endgroup\resethooks

The actual diagonal we return is a linear combination of faces. In
order to implement this, we use the notion of a \ensuremath{\Conid{\Conid{Data}.Map}}. This is an
associative array, indexing integer values representing coefficients
using the faces as indexing keys.

\begingroup\par\noindent\advance\leftskip\mathindent\(
\begin{pboxed}\SaveRestoreHook
\column{B}{@{}l@{}}
\column{3}{@{}l@{}}
\column{5}{@{}l@{}}
\column{E}{@{}l@{}}
\>[B]{}\mathbf{type}\;\Conid{LinearCombination}\;\Varid{vectors}\mathrel{=}\Conid{\Conid{Map}.Map}\;\Varid{vectors}\;\Conid{Int}{}\<[E]
\\[\blanklineskip]
\>[B]{}\Varid{showLinearCombination}\mathbin{::}\Conid{LinearCombination}\;\Conid{Face}\to \Conid{String}{}\<[E]
\\
\>[B]{}\Varid{showLinearCombination}\;\Varid{lc}\mathrel{=}\Varid{concatMap}\;\Varid{showSignFace}\;(\Varid{signFaceList}\;\Varid{lc}){}\<[E]
\\[\blanklineskip]
\>[B]{}\Varid{addSignFaces}\mathbin{::}[\mskip1.5mu \Conid{SignFace}\mskip1.5mu]\to \Conid{LinearCombination}\;\Conid{Face}{}\<[E]
\\
\>[B]{}\Varid{addSignFaces}\;[\mskip1.5mu \mskip1.5mu]\mathrel{=}\Varid{\Conid{Map}.empty}{}\<[E]
\\
\>[B]{}\Varid{addSignFaces}\;((\Varid{s},\Varid{p},\Varid{q})\mathbin{:}\Varid{as})\mathrel{=}\Varid{\Conid{Map}.insertWith}\;(\mathbin{+})\;(\Varid{p},\Varid{q})\;\Varid{s}\;(\Varid{addSignFaces}\;\Varid{as}){}\<[E]
\\[\blanklineskip]
\>[B]{}\Varid{signFaceList}\mathbin{::}\Conid{LinearCombination}\;\Conid{Face}\to [\mskip1.5mu \Conid{SignFace}\mskip1.5mu]{}\<[E]
\\
\>[B]{}\Varid{signFaceList}\;\Varid{lc}\mathrel{=}\Varid{map}\;(\lambda ((\Varid{p},\Varid{q}),\Varid{s})\to (\Varid{s},\Varid{p},\Varid{q}))\;(\Varid{\Conid{Map}.toList}\;\Varid{lc}){}\<[E]
\\[\blanklineskip]
\>[B]{}\Varid{permutahedronDiagonal}\mathbin{::}\Conid{Int}\to \Conid{LinearCombination}\;\Conid{Face}{}\<[E]
\\
\>[B]{}\Varid{permutahedronDiagonal}\;\Varid{n}\mathrel{=}(\Varid{addSignFaces}\mathbin{\circ}\Varid{nub})\;\Varid{deriveds}{}\<[E]
\\
\>[B]{}\hsindent{3}{}\<[3]
\>[3]{}\mathbf{where}{}\<[E]
\\
\>[3]{}\hsindent{2}{}\<[5]
\>[5]{}\Varid{deriveds}\mathrel{=}\Varid{concatMap}\;\Varid{derivedFaces}\;\Varid{primitiveFaces}{}\<[E]
\\
\>[3]{}\hsindent{2}{}\<[5]
\>[5]{}\Varid{primitiveFaces}\mathrel{=}\Varid{map}\;\Varid{buildFace}\;(\Varid{permutations}\;\Varid{n}){}\<[E]
\ColumnHook
\end{pboxed}
\)\par\noindent\endgroup\resethooks

\subsection{A diagonal on the associahedron}
\label{sec:diag-assoc}

Given the diagonal on permutahedra constructed above, we construct a
diagonal on the associahedra by applying the projection due to Andy
Tonks \cite{tonks97}. In practice, this projection eliminates, in
practice, all faces that contain partitions that are not \emph{derived
consecutive} in the sense that
$[\min\pi_j,\max\pi_j]\subset\bigcup_{i\leq j}\pi_i$. for all
$\pi_j\in\pi$.

This test is easily written, and thus the projection reduces to
extracting the summands that pass the test.

\begingroup\par\noindent\advance\leftskip\mathindent\(
\begin{pboxed}\SaveRestoreHook
\column{B}{@{}l@{}}
\column{3}{@{}l@{}}
\column{5}{@{}l@{}}
\column{7}{@{}l@{}}
\column{9}{@{}l@{}}
\column{13}{@{}l@{}}
\column{16}{@{}l@{}}
\column{33}{@{}l@{}}
\column{E}{@{}l@{}}
\>[B]{}\Varid{derivedConsecutive}\mathbin{::}\Conid{Partition}\to \Conid{Bool}{}\<[E]
\\
\>[B]{}\Varid{derivedConsecutive}\;\Varid{pi}\mathrel{=}\Varid{checkPartition}\;[\mskip1.5mu \mskip1.5mu]\;\Varid{pi}{}\<[E]
\\
\>[B]{}\hsindent{3}{}\<[3]
\>[3]{}\mathbf{where}{}\<[E]
\\
\>[3]{}\hsindent{2}{}\<[5]
\>[5]{}\Varid{checkPartition}\;\Varid{n}\;[\mskip1.5mu \mskip1.5mu]{}\<[33]
\>[33]{}\mathrel{=}\Conid{True}{}\<[E]
\\
\>[3]{}\hsindent{2}{}\<[5]
\>[5]{}\Varid{checkPartition}\;\Varid{n}\;(\Varid{pij}\mathbin{:}\Varid{pi'}){}\<[33]
\>[33]{}\mathrel{=}{}\<[E]
\\
\>[5]{}\hsindent{2}{}\<[7]
\>[7]{}\mathbf{if}\;{}\<[13]
\>[13]{}(\Varid{intersect}\;\Varid{n'}\;\Varid{range}\equiv \Varid{range}){}\<[E]
\\
\>[5]{}\hsindent{2}{}\<[7]
\>[7]{}\mathbf{then}\;{}\<[13]
\>[13]{}\Varid{checkPartition}\;\Varid{n'}\;\Varid{pi'}{}\<[E]
\\
\>[5]{}\hsindent{2}{}\<[7]
\>[7]{}\mathbf{else}\;{}\<[13]
\>[13]{}\Conid{False}{}\<[E]
\\
\>[5]{}\hsindent{2}{}\<[7]
\>[7]{}\mathbf{where}{}\<[E]
\\
\>[7]{}\hsindent{2}{}\<[9]
\>[9]{}\Varid{n'}{}\<[16]
\>[16]{}\mathrel{=}\Varid{sort}\;(\Varid{n}\plus \Varid{pij}){}\<[E]
\\
\>[7]{}\hsindent{2}{}\<[9]
\>[9]{}\Varid{range}{}\<[16]
\>[16]{}\mathrel{=}[\mskip1.5mu \Varid{minimum}\;\Varid{pij}\mathinner{\ldotp\ldotp}\Varid{maximum}\;\Varid{pij}\mskip1.5mu]{}\<[E]
\\[\blanklineskip]
\>[B]{}\Varid{associahedronDiagonal}\mathbin{::}\Conid{Int}\to \Conid{LinearCombination}\;\Conid{Face}{}\<[E]
\\
\>[B]{}\Varid{associahedronDiagonal}\;\Varid{n}\mathrel{=}{}\<[E]
\\
\>[B]{}\hsindent{3}{}\<[3]
\>[3]{}\Varid{\Conid{Map}.filterWithKey}\;\Varid{checkFace}\;(\Varid{permutahedronDiagonal}\;\Varid{n}){}\<[E]
\\
\>[B]{}\hsindent{3}{}\<[3]
\>[3]{}\mathbf{where}{}\<[E]
\\
\>[3]{}\hsindent{2}{}\<[5]
\>[5]{}\Varid{checkFace}\;(\Varid{f1},\Varid{f2})\;\Varid{s}\mathrel{=}{}\<[E]
\\
\>[5]{}\hsindent{2}{}\<[7]
\>[7]{}\Varid{derivedConsecutive}\;\Varid{f1}\mathrel{\wedge}\Varid{derivedConsecutive}\;(\Varid{reverse}\;\Varid{f2}){}\<[E]
\ColumnHook
\end{pboxed}
\)\par\noindent\endgroup\resethooks

\section{Installation, usage and example calculations}
\label{sec:inst-usage-example}

The way the paper is written, the source code of the program is the
source code of the paper. It can be downloaded with the paper source
code from \text{\tt arXiv:0707.4399}

In order to use this program, a Haskell interpreter or compiler would
have to be installed. I have used the Glasgow Haskell Compiler (GHC)
version 6.6.1 \cite{GHC} to develop this, but it should run on any
platform supporting the Haskell98 standard \cite{Haskell98}.

A typical worksession might look something like Example \ref{ex:usage}
\begin{example}[p]
\begin{tabbing}\tt
~\char36{}~ghci~SaneblidzeUmbleSigns\char46{}lhs\\
\tt ~~~~\char95{}\char95{}\char95{}~~~~~~~~~\char95{}\char95{}\char95{}~\char95{}\\
\tt ~~~\char47{}~\char95{}~\char92{}~\char47{}\char92{}~~\char47{}\char92{}\char47{}~\char95{}\char95{}\char40{}\char95{}\char41{}\\
\tt ~~\char47{}~\char47{}\char95{}\char92{}\char47{}\char47{}~\char47{}\char95{}\char47{}~\char47{}~\char47{}~~\char124{}~\char124{}~~~~~~GHC~Interactive\char44{}~version~6\char46{}6\char44{}~for~Haskell~98\char46{}\\
\tt ~\char47{}~\char47{}\char95{}\char92{}\char92{}\char47{}~\char95{}\char95{}~~\char47{}~\char47{}\char95{}\char95{}\char95{}\char124{}~\char124{}~~~~~~http\char58{}\char47{}\char47{}www\char46{}haskell\char46{}org\char47{}ghc\char47{}\\
\tt ~\char92{}\char95{}\char95{}\char95{}\char95{}\char47{}\char92{}\char47{}~\char47{}\char95{}\char47{}\char92{}\char95{}\char95{}\char95{}\char95{}\char47{}\char124{}\char95{}\char124{}~~~~~~Type~\char58{}\char63{}~for~help\char46{}\\
\tt ~\\
\tt ~Loading~package~base~\char46{}\char46{}\char46{}~linking~\char46{}\char46{}\char46{}~done\char46{}\\
\tt ~\char91{}1~of~1\char93{}~Compiling~SaneblidzeUmbleSigns~\char40{}~SaneblidzeUmbleSigns\char46{}lhs\char44{}~interpreted~\char41{}\\
\tt ~Ok\char44{}~modules~loaded\char58{}~SaneblidzeUmbleSigns\char46{}\\
\tt ~\char42{}SaneblidzeUmbleSigns\char62{}~Map\char46{}size~\char46{}~permutahedronDiagonal~\char36{}~4\\
\tt ~50\\
\tt ~\char42{}SaneblidzeUmbleSigns\char62{}~Map\char46{}size~\char46{}~associahedronDiagonal~\char36{}~4\\
\tt ~22\\
\tt ~\char42{}SaneblidzeUmbleSigns\char62{}~map~\char40{}Map\char46{}size~\char46{}~permutahedronDiagonal\char41{}~\char91{}1\char46{}\char46{}6\char93{}\\
\tt ~\char91{}1\char44{}2\char44{}8\char44{}50\char44{}432\char44{}4802\char93{}\\
\tt ~\char42{}SaneblidzeUmbleSigns\char62{}~map~\char40{}Map\char46{}size~\char46{}~associahedronDiagonal\char41{}~\char91{}1\char46{}\char46{}6\char93{}\\
\tt ~\char91{}1\char44{}2\char44{}6\char44{}22\char44{}91\char44{}408\char93{}\\
\tt ~\char42{}SaneblidzeUmbleSigns\char62{}~putStr~\char46{}~unlines~\char46{}~map~showSignFace~\char46{}~\\
\tt ~~~~~~~~~~signFaceList~\char46{}~permutahedronDiagonal~\char36{}~3\\
\tt ~\char43{}1\char124{}2\char124{}3x1\char44{}2\char44{}3\\
\tt ~\char45{}1\char124{}2\char44{}3x3\char124{}1\char44{}2\\
\tt ~\char45{}1\char124{}2\char44{}3x1\char44{}3\char124{}2\\
\tt ~\char43{}1\char44{}2\char124{}3x2\char44{}3\char124{}1\\
\tt ~\char43{}1\char44{}2\char124{}3x2\char124{}1\char44{}3\\
\tt ~\char43{}1\char44{}2\char44{}3x3\char124{}2\char124{}1\\
\tt ~\char45{}1\char44{}3\char124{}2x3\char124{}1\char44{}2\\
\tt ~\char43{}2\char124{}1\char44{}3x2\char44{}3\char124{}1\\
\tt ~\char42{}SaneblidzeUmbleSigns\char62{}~putStr~\char46{}~unlines~\char46{}~map~showSignFace~\char46{}~\\
\tt ~~~~~~~~~~signFaceList~\char46{}~associahedronDiagonal~\char36{}~3\\
\tt ~\char43{}1\char124{}2\char124{}3x1\char44{}2\char44{}3\\
\tt ~\char45{}1\char124{}2\char44{}3x3\char124{}1\char44{}2\\
\tt ~\char43{}1\char44{}2\char124{}3x2\char44{}3\char124{}1\\
\tt ~\char43{}1\char44{}2\char124{}3x2\char124{}1\char44{}3\\
\tt ~\char43{}1\char44{}2\char44{}3x3\char124{}2\char124{}1\\
\tt ~\char43{}2\char124{}1\char44{}3x2\char44{}3\char124{}1\\
\tt ~\char42{}SaneblidzeUmbleSigns\char62{}~\char58{}q\\
\tt ~Leaving~GHCi\char46{}
\end{tabbing}
\caption{Some calculations with the \text{\tt SaneblidzeUmble\char46{}lhs}
  implementation}
\label{ex:usage}
\end{example}

\subsection{Performance data}
\label{sec:performance-data}

We have tested the code and its performance by calculating,
subsequently, the diagonals on $P_1,\dots,P_7$. The results of our
calculations as well as some complexity measurements can be found in
Table \ref{tbl:performance}. The tests were performed on a double Dual
Core AMD Opteron(tm) Processor 270 with 16G RAM running OpenSuSE 10.2
with standard linux kernel version 2.6.18.

Given the garbage collection that GHC uses, there is a difference to
be observed between the total amount of memory ever allocated, and the
maximal amount of memory allocated at a single point in time. The
measurements will state both. 

\begin{table}
  \centering
 \begin{tabular}{r|r|r@{.}l|r@{.}l|r@{.}l}
   $n$ & $|\Delta_{P_n}(e^n)|$ & 
\multicolumn{2}{c}{Execution time} &
\multicolumn{2}{c}{Total allocation} &
\multicolumn{2}{c}{Peak allocation} \\
   \hline
1 &   $1$                  & $<$0&005s  &  45&813k  &  28&617k   \\
2 &   $2$                  & $<$0&005s  &  53&438k  &  28&617k   \\
3 &   $8$                  & $<$0&005s  &  95&648k  &  28&617k   \\
4 &   $50$                 & $<$0&005s  & 402&500k  &  28&617k   \\
5 &   $432$                &    0&02s   &   3&987M  &  56&695k   \\
6 &   $4\,802$             &    1&63s   &  45&687M  &   1&631M   \\
7 &   $65\,536$            &  399&97s   &   1&000G  &  22&198M   \\
8 &   $1\,062\,882$        &93$\,$965&64s& 39&205G  & 342&704M   \\
9 &   $20\,000\,000$       & 
\multicolumn{2}{c}{$>$400h} &
\multicolumn{2}{c}{N/A} &
$>$3&000G
 \end{tabular}
  \caption{Performance and calculations}
  \label{tbl:performance}
\end{table}

%\clearpage

\section{Haskell notation}
\label{sec:haskell-notation}

The paper uses Haskell code for the algorithm discussion, and for the
benefit of the reader, we shall here discuss the notation we used to
provide a dictionary and help elucidate the code.

\subsection{Types and definitions}
\label{sec:types-definitions}

Haskell is a statically typed language, which means that any object
has a type, and can only be used in an expression if the type matches
all functions involved. Types can be constructed from a family of
primitive types, of which we here see \ensuremath{\Conid{Int}} and \ensuremath{\Conid{Bool}} occuring. 

\ensuremath{\Conid{Int}} is a machine word type, admitting integers up to about $2^{32}$
before overflowing. Since the code runs into performance problems
before the integers we use reach $10$, this is not viewed as a problem
at this point.

\ensuremath{\Conid{Bool}} is a truth value type, with the two valid elements \ensuremath{\Conid{True}} and
\ensuremath{\Conid{False}}. 

Using these types, we are then able to construct more complex types.
The most important derived type templates we use are the list types,
the \ensuremath{\Conid{Maybe}} type, maps and the pair type. We shall deal with each of
these constructions in a section of their own.

All functions we declare have their type signature given
explicitly. This is done by a declaration on the form
\begingroup\par\noindent\advance\leftskip\mathindent\(
\begin{pboxed}\SaveRestoreHook
\column{B}{@{}l@{}}
\column{3}{@{}l@{}}
\column{E}{@{}l@{}}
\>[3]{}\Varid{functionName}\mathbin{::}\Varid{firstArgType}\to \Varid{secondArgType}\to \mathbin{...}\to \Varid{returnType}{}\<[E]
\ColumnHook
\end{pboxed}
\)\par\noindent\endgroup\resethooks

The function body then is declared by giving names to the function
arguments, and declaring what happens with them. While doing this,
there are two specific constructions with special meaning. First,
there is the \ensuremath{\anonymous } variable name, which means that that particular input
variable is not named in the coming declaration. This is, for instance
used, for cases where certain of the argument are not needed, for
instance in the \ensuremath{\Varid{monotonicSequence}} definition, it doesn't matter what
comparison function is used if the sequence is empty or a
singleton. Thus the definitions 
\begingroup\par\noindent\advance\leftskip\mathindent\(
\begin{pboxed}\SaveRestoreHook
\column{B}{@{}l@{}}
\column{26}{@{}l@{}}
\column{E}{@{}l@{}}
\>[B]{}\Varid{monotonicSequence}\;\anonymous \;[\mskip1.5mu \mskip1.5mu]{}\<[26]
\>[26]{}\mathrel{=}[\mskip1.5mu \mskip1.5mu]{}\<[E]
\\
\>[B]{}\Varid{monotonicSequence}\;\anonymous \;[\mskip1.5mu \Varid{x}\mskip1.5mu]{}\<[26]
\>[26]{}\mathrel{=}[\mskip1.5mu [\mskip1.5mu \Varid{x}\mskip1.5mu]\mskip1.5mu]{}\<[E]
\ColumnHook
\end{pboxed}
\)\par\noindent\endgroup\resethooks

The other special character in argument definitions is \ensuremath{\mathord{@}}. It can be
used to allocate a name to an entire complex construct, while
allocating names to the entries as well. This is done in the output
function \ensuremath{\Varid{showSignFace}}, as follows
\begingroup\par\noindent\advance\leftskip\mathindent\(
\begin{pboxed}\SaveRestoreHook
\column{B}{@{}l@{}}
\column{E}{@{}l@{}}
\>[B]{}\Varid{showSignFace}\;\Varid{f}\mathord{@}(\Varid{s},\anonymous ,\anonymous )\mathrel{=}\mathbin{...}{}\<[E]
\ColumnHook
\end{pboxed}
\)\par\noindent\endgroup\resethooks
whereby here \ensuremath{\Varid{f}} ends up refering to the entire \ensuremath{\Conid{SignFace}} it acts on,
and \ensuremath{\Varid{s}} additionally refers to the sign of that face. The two \ensuremath{\anonymous }
indicate that we do not wish to have names allocated to the partitions
of the face.

Finally, we have the ability to construct type synonyms. By a
declaration like
\begingroup\par\noindent\advance\leftskip\mathindent\(
\begin{pboxed}\SaveRestoreHook
\column{B}{@{}l@{}}
\column{3}{@{}l@{}}
\column{E}{@{}l@{}}
\>[3]{}\mathbf{type}\;\Conid{Newtype}\mathrel{=}\Conid{Oldtype}{}\<[E]
\ColumnHook
\end{pboxed}
\)\par\noindent\endgroup\resethooks
we can henceforth use \ensuremath{\Conid{Newtype}} as a type in its own right, which can
be used anywhere where \ensuremath{\Conid{Oldtype}} could be used. We use this to make
readable code, thus we can write, for instance
\begingroup\par\noindent\advance\leftskip\mathindent\(
\begin{pboxed}\SaveRestoreHook
\column{B}{@{}l@{}}
\column{3}{@{}l@{}}
\column{19}{@{}l@{}}
\column{E}{@{}l@{}}
\>[3]{}\mathbf{type}\;\Conid{Sequence}{}\<[19]
\>[19]{}\mathrel{=}[\mskip1.5mu \Conid{Int}\mskip1.5mu]{}\<[E]
\\
\>[3]{}\mathbf{type}\;\Conid{Partition}{}\<[19]
\>[19]{}\mathrel{=}[\mskip1.5mu \Conid{Sequence}\mskip1.5mu]{}\<[E]
\\
\>[3]{}\mathbf{type}\;\Conid{Face}{}\<[19]
\>[19]{}\mathrel{=}(\Conid{Partition},\Conid{Partition}){}\<[E]
\\
\>[3]{}\Varid{derivedFaces}\mathbin{::}\Conid{Face}\to [\mskip1.5mu \Conid{Face}\mskip1.5mu]{}\<[E]
\ColumnHook
\end{pboxed}
\)\par\noindent\endgroup\resethooks
instead of the more unwieldy und much less readable
\begingroup\par\noindent\advance\leftskip\mathindent\(
\begin{pboxed}\SaveRestoreHook
\column{B}{@{}l@{}}
\column{3}{@{}l@{}}
\column{E}{@{}l@{}}
\>[3]{}\Varid{derivedFaces}\mathbin{::}([\mskip1.5mu [\mskip1.5mu \Conid{Int}\mskip1.5mu]\mskip1.5mu],[\mskip1.5mu [\mskip1.5mu \Conid{Int}\mskip1.5mu]\mskip1.5mu])\to [\mskip1.5mu ([\mskip1.5mu [\mskip1.5mu \Conid{Int}\mskip1.5mu]\mskip1.5mu],[\mskip1.5mu [\mskip1.5mu \Conid{Int}\mskip1.5mu]\mskip1.5mu])\mskip1.5mu]{}\<[E]
\ColumnHook
\end{pboxed}
\)\par\noindent\endgroup\resethooks

\subsection{Lists and their manipulation}
\label{sec:lists-their-manip}

A list type takes the form, on the type level, of \ensuremath{[\mskip1.5mu \Conid{Type}\mskip1.5mu]}. It
consists of a randomly accessable sequence of entries, all of which
has to have type \ensuremath{\Conid{Type}}. Lists can be given literally, by enclosing
the comma-separated list of entries in square brackets, like
\ensuremath{[\mskip1.5mu \mathrm{2},\mathrm{3},\mathrm{4}\mskip1.5mu]}.  An element can be prepended to a list using \ensuremath{\mathbin{:}}, for
example, \ensuremath{\mathrm{1}\mathbin{:}[\mskip1.5mu \mathrm{2},\mathrm{3},\mathrm{4}\mskip1.5mu]\mathrel{=}[\mskip1.5mu \mathrm{1},\mathrm{2},\mathrm{3},\mathrm{4}\mskip1.5mu]}. For the reader with Lisp experience,
\ensuremath{\mathbin{:}} is very much like \ensuremath{\Varid{cons}}. The list \ensuremath{[\mskip1.5mu \mskip1.5mu]} is a special case -- the
empty list -- and corresponds in Lisp terms to \ensuremath{\Conid{Nil}}.  The integers
between \ensuremath{\Varid{a}} and \ensuremath{\Varid{b}} can be listed by the construction \ensuremath{[\mskip1.5mu \Varid{a}\mathinner{\ldotp\ldotp}\Varid{b}\mskip1.5mu]}.

Note that strings are really lists of type \ensuremath{[\mskip1.5mu \Conid{Char}\mskip1.5mu]}, so that \ensuremath{\text{\tt \char34 abc\char34}\mathrel{=}[\mskip1.5mu \text{\tt 'a'},\text{\tt 'b'},\text{\tt 'c'}\mskip1.5mu]}. This allows us to handle strings with the kind of
toolbox we build for list manipulation as well.

When defining functions, we make large use of pattern matching on the
input parameters. Thus, by stating several special cases and then a
general case, we can achieve very compact and readable
definitions. For this, the most important constructions are
\begin{description}
\item[\ensuremath{[\mskip1.5mu \mskip1.5mu]}] The empty list.
\item[\ensuremath{\anonymous }] Any input parameter at all -- we simply do not care for the
  value of this parameter, and will not use it in the definition.
\item[\ensuremath{[\mskip1.5mu \Varid{a}\mskip1.5mu]}] A singleton list, containing the element \ensuremath{\Varid{a}}.
\item[\ensuremath{(\Varid{a}\mathbin{:}\Varid{as})}] The element \ensuremath{\Varid{a}} followed by the rest of the list
  \ensuremath{\Varid{as}}. This can be expanded, so that \ensuremath{(\Varid{a}\mathbin{:}\Varid{b}\mathbin{:}\Varid{c}\mathbin{:}\Varid{etc})} picks out the
  first three elements of a list with at least three elements, and
  stores them in the variables \ensuremath{\Varid{a}}, \ensuremath{\Varid{b}} and \ensuremath{\Varid{c}} for the duration of
  the function definition.
\end{description}

Apart from the pattern matching, there is a wealth of functions
available for handling lists. The ones used in the code here are
\begin{description}
\item[\ensuremath{\Varid{map}}] applies a given function to all elements of a given list.
\item[\ensuremath{\Varid{foldl'}}] takes a binary function $*$, a first argument $a$ and a list of arguments
  $a_n$, and returns $(\dots(((a*a_1)*a_2)*a_3)\dots)*a_n$.
\item[\ensuremath{\Varid{reverse}}] reverses a finite list.
\item[\ensuremath{\plus }] concatenates two lists.
\item[\ensuremath{\Varid{intersperse}}] takes an element and inserts it between all
  elements of a list.
\item[\ensuremath{\Varid{concat}}] flattens a list of lists one level.
\item[\ensuremath{\Varid{concatMap}}] first applies a function to all elements of a list,
  and then flattens the result.
\item[\ensuremath{\Varid{unlines}}] intersperses the newline character, and concatenates
  the result to give a string representing each string in the list it
  is applied to on a new line on its own.
\item[\ensuremath{\Varid{head}}] returns the first element of a list.
\item[\ensuremath{\Varid{length}}] returns the number of elements of a list.
\item[\ensuremath{\Varid{intersect}}] returns the list of elements that occur in both the
  argument lists.
\item[\ensuremath{\Varid{list}\mathbin{!!}\Varid{i}}] returns the element of index \ensuremath{\Varid{i}}. Thus \ensuremath{\Varid{head}\;\Varid{list}\equiv \Varid{list}\mathbin{!!}\mathrm{0}}.
\item[\ensuremath{\Varid{minimum}}] returns the minimal element of a list of comparable
  elements.
\item[\ensuremath{\Varid{maximum}}] returns the maximal element of a list of comparable
  elements.
\item[\ensuremath{\Varid{findIndices}}] returns a list of all \ensuremath{\Varid{i}} such that the function
  given evaluates to \ensuremath{\Conid{True}} for the element of index \ensuremath{\Varid{i}}, and \ensuremath{\Conid{False}}
  for all other elements of the list.
\item[\ensuremath{\Varid{findIndex}}] returns \ensuremath{\Conid{Just}\;\Varid{i}} if \ensuremath{\Varid{i}} is the first element of the
  result from \ensuremath{\Varid{findIndices}}, and \ensuremath{\Conid{Nothing}} if \ensuremath{\Varid{findIndices}} would
  return an empty list.
\item[\ensuremath{\Varid{null}}] returns \ensuremath{\Conid{True}} if the argument is \ensuremath{[\mskip1.5mu \mskip1.5mu]} and \ensuremath{\Conid{False}}
  otherwise.
\item[\ensuremath{\Varid{list}\mathbin{\char92 \char92 }\Varid{list'}}] removes all elements in \ensuremath{\Varid{list'}} from the list
  \ensuremath{\Varid{list}}.
\item[\ensuremath{\Varid{take}\;\Varid{j}\;\Varid{list}}] returns the list of the first \ensuremath{\Varid{j}} elements of 
  \ensuremath{\Varid{list}}.
\item[\ensuremath{\Varid{drop}\;\Varid{j}\;\Varid{list}}] returns the list of the elements following the
  first \ensuremath{\Varid{j}} elements of \ensuremath{\Varid{list}}.
\item[\ensuremath{\Varid{filter}}] returns only those elements for which the function
  given evaluates to \ensuremath{\Conid{True}}.
\item[\ensuremath{\Varid{tails}\;\Varid{list}}] is equivalent to \ensuremath{\Varid{map}\;(\lambda \Varid{j}\to \Varid{drop}\;\Varid{j}\;\Varid{list})\;[\mskip1.5mu \mathrm{1}\mathinner{\ldotp\ldotp}\Varid{length}\;\Varid{list}\mskip1.5mu]}.
\item[\ensuremath{\Varid{nub}}] removes duplicates from a list.
\item[\ensuremath{\Varid{sort}}] sorts the list.
\item[\ensuremath{\Varid{elem}}] checks whether an element appears in a list.
\item[\ensuremath{\Varid{sum}}] sums the elements of a list.
\end{description}

\subsection{Maybe -- graceful error handling}
\label{sec:maybe-graceful-error}

It is quite possible that a calculation can, or may, be expected to
fail on some input. The canonical way to handle this in Haskell is to
assign a type to the results of the calculation that handles error
cases. Unless an error description is necessary, the \ensuremath{\Conid{Maybe}} type
encapsulation is what gets used. The \ensuremath{\Conid{Maybe}\;\Conid{Mytype}} has two different
kinds of elements -- either \ensuremath{\Conid{Nothing}} or \ensuremath{\Conid{Just}\;\Varid{aValue}}. The return
value of \ensuremath{\Conid{Nothing}}, thus, is used to signify a failed calculation.

In order to handle output from this kind of calculation, there are a
few methods available.
\begin{description}
\item[Pattern matching] with the patterns \ensuremath{\Conid{Nothing}} or \ensuremath{\Conid{Just}\;\Varid{a}} can be
  used to separate into cases depending on the results of the calculation.
\item[\ensuremath{\Varid{isNothing}}] returns \ensuremath{\Conid{True}} if the object considered is
  \ensuremath{\Conid{Nothing}}.
\item[\ensuremath{\Varid{isJust}}] returns \ensuremath{\Conid{True}} if the object considered is not \ensuremath{\Conid{Nothing}}.
\item[\ensuremath{\Varid{fromJust}}] takes \ensuremath{\Conid{Just}\;\Varid{a}} and returns \ensuremath{\Varid{a}}. If used on \ensuremath{\Conid{Nothing}}
  causes an exception to be thrown, thus aborting execution
  entirely. To be used only when it is impossible to receive a \ensuremath{\Conid{Nothing}}.
\item[\ensuremath{\Varid{mapMaybe}}] maps a function returning \ensuremath{\Conid{Maybe}\;\Conid{Sometype}} over a
  list, and returns a list of only the results that were encapsulated
  in \ensuremath{\Conid{Just}}.
\end{description}

\subsection{Pairs of elements}
\label{sec:pairs-elements}

The expression \ensuremath{(\Varid{a},\Varid{b})} gives an element of type \ensuremath{(\Conid{AType},\Conid{BType})}. This
is a way to construct static length tuples of elements that need not
be of the same type -- as different from the lists, where all elements
must be of the same type. The functions \ensuremath{\Varid{fst}} and \ensuremath{\Varid{snd}} allow access
to \ensuremath{\Varid{a}} and \ensuremath{\Varid{b}} respectively. Pattern matching works in the expected
way -- \ensuremath{(\Varid{a},\Varid{b})} assigns the name \ensuremath{\Varid{a}} to the first element and \ensuremath{\Varid{b}} to
the second in a pattern match.

\subsection{\ensuremath{\Conid{\Conid{Data}.Map}} -- associative arrays}
\label{sec:data.m-assoc-arrays}

An associative array is a way to store data with other indexing sets
than just the non-negative integers. Thus, we can use (almost) any
indexing type in an associative array. The application we use it for
here is for the type carrying linear combinations of faces. Here, a
linear combination is an allocation from the set of occuring faces to
the coefficients for each face.

In order to interact with these constructs, there are a number of
operations available -- many of which are not used in the code, but
will prove useful to any user who wants to extract information from a
calculation.

\begin{description}
\item[\ensuremath{\Varid{signFaceList}}] converts back from the linear combinations to a
  list of signed faces.
\item[\ensuremath{\Varid{showLinearCombination}}] produces a string representing the
  entire linear combination.
\item[\ensuremath{\Varid{\Conid{Map}.size}}] gives the number of terms in a linear combination,
  though not necessarily excluding terms with coefficient 0.
\item[\ensuremath{\Varid{\Conid{Map}.filter}}] removes terms for which the value held evaluates to
  \ensuremath{\Conid{False}} in a given function. Thus \ensuremath{\Varid{\Conid{Map}.filter}\;(\not\equiv \mathrm{0})} can be used to
  remove all trivial terms from the results.
\item[\ensuremath{(\Varid{permutahedronDiagonal}\;\Varid{n})\mathbin{!}\Varid{someFace}}] returns the sign of
  \ensuremath{\Varid{someFace}}.
\item[\ensuremath{\Varid{\Conid{Map}.elems}}] returns all the coefficients in the diagonal.
\item[\ensuremath{\Varid{\Conid{Map}.keys}}] returns all the faces in the diagonal.
\end{description}

\subsection{Logical and arithmetic operations}
\label{sec:logic-arithm-oper}

Logic has much of the expected operators, thus we have
\begin{description}
\item[\ensuremath{\mathrel{\wedge}}] for logical and,
\item[\ensuremath{\mid }] for logical or,
\item[\ensuremath{\neg }] for logical negation,
\item[\ensuremath{\equiv }] for equality testing,
\item[\ensuremath{\not\equiv }] for inequality testing,
\item[\ensuremath{\geq }, \ensuremath{\mathbin{>}}, \ensuremath{\leq } and \ensuremath{\mathbin{<}}] for orderings,
\item[\ensuremath{\mathbin{+}}] for addition,
\item[\ensuremath{\mathbin{*}}] for multiplication,
\item[\ensuremath{\mathbin{\uparrow}}] for exponentiation and
\item[\ensuremath{\mathbin{-}}] for subtraction or negation, depending on context.
\end{description}

\subsection{Loose ends}
\label{sec:loose-ends}

Finally, we mention several purely syntactical subtleties.
The keyword \ensuremath{\mathbf{where}} can be used to postpone definition of symbols used
until an expression is completely given. Thus, anything listed after
the \ensuremath{\mathbf{where}} becomes definitions that scope over the expression \ensuremath{\mathbf{where}}
follows.

The \ensuremath{\mathbf{if}}-\ensuremath{\mathbf{then}}-\ensuremath{\mathbf{else}} constructs give little surprise, but are a language
construct of their own, specializing expressions on the form
\begingroup\par\noindent\advance\leftskip\mathindent\(
\begin{pboxed}\SaveRestoreHook
\column{B}{@{}l@{}}
\column{3}{@{}l@{}}
\column{5}{@{}l@{}}
\column{E}{@{}l@{}}
\>[3]{}\mathbf{case}\;\Varid{something}\;\mathbf{of}{}\<[E]
\\
\>[3]{}\hsindent{2}{}\<[5]
\>[5]{}\Varid{value1}\to \Varid{result1}{}\<[E]
\\
\>[3]{}\hsindent{2}{}\<[5]
\>[5]{}\Varid{value2}\to \Varid{result2}{}\<[E]
\ColumnHook
\end{pboxed}
\)\par\noindent\endgroup\resethooks
where the choice of result is made over a larger array of alternatives
than only \ensuremath{\Conid{True}} and \ensuremath{\Conid{False}}.

The special operators \ensuremath{\mathbin{\circ}} and \ensuremath{\mathbin{\$}} %$
handle function composition. The operation \ensuremath{\mathbin{\circ}} works entirely as
expected, with \ensuremath{\Varid{f}\;(\Varid{g}\;(\Varid{h}\;\Varid{x}))\mathrel{=}(\Varid{f}\mathbin{\circ}\Varid{g}\mathbin{\circ}\Varid{h})\;\Varid{x}}. \ensuremath{\mathbin{\$}} %$
on the other hand works more as function application, and has a kind
of backwards associativity, so that instead of writing \ensuremath{\Varid{f}\;(\Varid{someLargeExpression})} we can content ourselves by writing \ensuremath{\Varid{f}\mathbin{\$}\Varid{someLargeExpression}}. %$

Since, in functional programming, functions are just as good arguments
to other functions as any other sort of value, there is a need for
some way to construct an anonymous, lightweight, throwaway function to
be inserted in some particular context or other. This is for
historical reasons done with the \ensuremath{\lambda } construction. By an expression of
the form \ensuremath{\lambda \Varid{x}\to \Varid{someExpression}\;\Varid{x}} we declare a function that takes an
argument \ensuremath{\Varid{x}} and returns \ensuremath{\Varid{someExpression}\;\Varid{x}}. 

\section{Acknowledgements}
\label{sec:acknowledgements}

This paper would not have been possible without all the fruitful
conversations held with Ron Umble. The author is thankful for these
and the generous help given in preparing this preprint.

The author further acknowledges that the ideas would have remained
largely unpublished without the urging by Jim Stasheff to write it up.

The Haskell section has enjoyed the benefits of several experienced
Haskell programmers, most notably Shae Erisson. 

The author would like to acknowledge his supervisor David J. Green,
who has remained supportive of the author's ideas and efforts.
\bibliography{library} 
\bibliographystyle{alpha} 
\end{document}